\documentclass[1p]{elsarticle}
\biboptions{sort&compress}
\usepackage{amsmath,amssymb,bm}
\usepackage{amsthm}
\usepackage{color}
\usepackage{hyperref}
\hypersetup{
    colorlinks=true,
    linkcolor=blue,     citecolor=red,     urlcolor=blue  }
\usepackage[colorinlistoftodos]{todonotes}
\usepackage{mathptmx}   
\usepackage[weather,alpine,misc,geometry]{ifsym}
\usepackage{amsfonts}
\usepackage{graphicx}%
\usepackage[capitalise,nameinlink]{cleveref}
\crefformat{equation}{(#2#1#3)}
\crefmultiformat{equation}{(#2#1#3)}%
{ and~(#2#1#3)}{, (#2#1#3)}{ and~(#2#1#3)}
\crefrangeformat{equation}{(#3#1#4)--(#5#2#6)}
\crefrangemultiformat{equation}{(#3#1#4)--(#5#2#6)}%
{ and~(#3#1#4)--(#5#2#6)}{, (#3#1#4)--(#5#2#6)}{ and~(#3#1#4)--(#5#2#6)}

\newcommand{\al}{\alpha}
\newcommand{\be}{\beta}
\newcommand{\ga}{\gamma}
\newcommand{\la}{\lambda}
\newcommand{\de}{\delta}
\newcommand{\eps}{\varepsilon}
\newcommand{\bx}{\bar x}
\newcommand{\by}{\bar y}

\newcommand {\R} {\mathbb R}
\newcommand {\N} {\mathbb N}

\newcommand {\B} {\mathbb B}

\newcommand{\dist}{{\rm dist}\,}
\newcommand {\gph} {{\rm gph}\,}
\newcommand {\dom} {{\rm dom}\,}
\newcommand {\epi} {{\rm epi}\,}

\newcommand {\cl} {{\rm cl}\,}
\newcommand {\co} {{\rm co}\,}

\newcommand {\sd} {\partial}
\newcommand {\Int} {{\rm int}\,}

\newcommand{\toto}{\rightrightarrows}
\newcommand{\folgt}{$ \Rightarrow\ $}

\def\nbh{neighbourhood}
\def\es{\emptyset}
\def\lsc{lower semicontinuous}

\def\LHS{left-hand side}
\def\RHS{right-hand side}

\def\EVP{Ekeland variational principle}
\def\Fr{Fr\'echet}

\newcommand{\norm}[1]{\left\Vert#1\right\Vert}

\newcommand{\ang}[1]{\left\langle #1 \right\rangle}
\newcommand{\qdtx}[1]{\quad\mbox{#1}\quad}
\newcommand{\AND}{\quad\mbox{and}\quad}
\newcounter{mycount}
\def\cnta{\setcounter{mycount}{\value{enumi}}}
\def\cntb{\setcounter{enumi}{\value{mycount}}}


\newtheorem{theorem}{Theorem}[section]
\newtheorem{proposition}[theorem]{Proposition}
\newtheorem{lemma}[theorem]{Lemma}
\newtheorem{corollary}[theorem]{Corollary}
\newtheorem{corollary.pr}{Corollary}

\theoremstyle{definition}
\newtheorem{definition}[theorem]{Definition}
\theoremstyle{remark}
\newtheorem{remark}[theorem]{Remark}
\newtheorem{example}[theorem]{Example}
\theoremstyle{plain}

\newcommand{\EI}{EI^{\textup{cl}}}
\usepackage{enumitem}
\setlist[enumerate,1]{label={\rm(\roman*)}}
\setlist[enumerate,2]{label={\rm(\alph*)}}

\newcommand{\weakly}{\rightharpoonup}

\newcommand{\blambda}{{\boldsymbol{\lambda}}}

\renewcommand{\es}{\emptyset}

\newcommand{\sdc}{{\partial}^{\textup{C}}}
\newcommand{\sdf}{{\partial}}
\newcommand{\sdm}{{\overline{\partial}}}

\DeclareMathAlphabet{\mathpzc}{OT1}{pzc}{m}{it}
\newcommand\oo{\mathpzc{o}}

\numberwithin{equation}{section}

\begin{document}

\begin{frontmatter}
\title{%
Fuzzy multiplier, sum and intersection rules in non-Lipschitzian settings:
Decoupling approach revisited\tnoteref{dedic}
}%
\tnotetext[dedic]{
Dedicated to Alexander Ioffe on the occasion of his 85$^{th}$ birthday}

	
\author[1]{Mari\'{a}n Fabian}
\affiliation[1]{organization={								Czech Academy of Sciences, Institute of Mathematics},
addressline={{Z}itna 25},
city={Prague 1},
postcode={115 67},
country={Czech Republic, ORCID: 0000-0003-3031-0862}}

\ead{fabian@math.cas.cz}
	
\author[2]{Alexander Y.\ Kruger\corref{cor1}}
\cortext[cor1]{Corresponding author}
\affiliation[2]{organization={								Optimization Research Group,
Faculty of Mathematics and Statistics,
Ton Duc Thang University},
addressline={19 Nguyen Huu Tho St, Tan Phong Ward, Dist. 7},
city={Ho Chi Minh City},
country={Vietnam, ORCID: 0000-0002-7861-7380}}

\ead{alexanderkruger@tdtu.edu.vn}

\author[3]{Patrick Mehlitz}
\affiliation[3]{
	organization	=	{Faculty of Mathematics, University of Duisburg-Essen},
	city			=	{Essen},
	postcode		=	{45127},
	country			=	{Germany, ORCID: 0000-0002-9355-850X}}
	\ead{patrick.mehlitz@uni-due.de}

\begin{abstract}
We revisit the \emph{decoupling approach} widely used (often intuitively) in nonlinear analysis and optimization
and initially formalized about a quarter of a century ago by Borwein \& Zhu, Borwein \& Ioffe and Lassonde.
It allows one to streamline proofs of necessary optimality conditions and calculus relations,
unify and simplify the respective statements, clarify and in many cases weaken the assumptions.
In this paper we study weaker concepts of \emph{quasiuniform infimum},
\emph{quasiuniform lower semicontinuity} and
\emph{quasiuniform minimum},
putting them into the context of the general theory developed by the aforementioned authors.
Along the way, we unify the terminology and notation and fill in some gaps in the general theory.
We establish rather general primal and dual necessary conditions characterizing quasiuniform
$\eps$-minima of the sum of two functions.
The obtained fuzzy multiplier rules are formulated in general Banach spaces in terms of Clarke subdifferentials
and in Asplund spaces in terms of \Fr\ subdifferentials.
The mentioned fuzzy multiplier rules naturally lead to certain fuzzy subdifferential calculus results.
An application from sparse optimal control illustrates applicability of the obtained findings.
\end{abstract}

\end{frontmatter}


\setlist{nosep}

\section{Introduction}\label{sec:introduction}

When dealing with problems involving several component functions or sets, e.g.,
proving necessary optimality conditions in metric or normed spaces
or establishing subdifferential/normal cone calculus relations in normed spaces,
it is common to consider extensions of the problems allowing the components
to depend on or involve individual variables while ensuring that these individual variables are not too far apart.
This \emph{decoupling method} (the term coined by Borwein and Zhu \cite{BorZhu05}) allows one
to express the resulting conditions in terms of subdifferentials of the individual functions
and/or normal cones to individual sets, or appropriate primal space tools.

For instance, when dealing with the problem of minimizing the sum of two extended-real-valued functions $\varphi_1$ and $\varphi_2$,
one often replaces the function of a single variable $x\mapsto(\varphi_1+\varphi_2)(x)$
with the function of two variables (\emph{decoupled sum} \cite{Las01})
$(x_1,x_2)\mapsto\varphi_1(x_1)+\varphi_2(x_2)$ while forcing the distance $d(x_1,x_2)$ to be small.
The latter is often done by adding a penalty term (or an increasing sequence of terms) containing $d(x_1,x_2)$.
Here and throughout the paper, for brevity, we restrict ourselves to the case of two functions.
The definitions and statements can be easily extended to an arbitrary finite number of functions.

The decoupling approach has been intuitively used in numerous publications for decades.
As claimed in \cite[Section~6.1.4]{BorZhu05}, all basic subdifferential rules in Banach spaces are
different facets of a variational principle in conjunction with a decoupling method.
The basics of the decoupling method were formalized in a series of publications
by Borwein and Zhu \cite{BorweinZhu96,BorZhu05}, Borwein and Ioffe \cite{BorweinIoffe1996}, Lassonde \cite{Las01} and Penot \cite{Pen13}.
With regards to the mentioned above minimization problem, the following \emph{uniform infimum} \cite{Las01}
\begin{align}
\label{La}
\Lambda_U(\varphi_1,\varphi_2):=& \liminf_{\substack{d(x_1,x_2)\to0,\,\dist(x_1,U)\to0}} (\varphi_1(x_1)+\varphi_2(x_2))
\end{align}
of $(\varphi_1,\varphi_2)$ over (or around) $U$ plays a key role.
Here $U$ is a given set.
It can represent a set of constraints or a \nbh\ of a given point.
Observe that, thanks to $d(x_1,x_2)\to0$, condition $\dist(x_1,U)\to0$ is equivalent to $\dist(x_2,U)\to0$;
hence definition \eqref{La} gives no advantage to the variable $x_1$.
The quantity $\Lambda_U(\varphi_1,\varphi_2)$ from \eqref{La} is referred to in \cite{BorZhu05} as \emph{decoupled infimum}
and in \cite{Pen13} as \emph{stabilized infimum}.
As pointed out in \cite{Las01,BorZhu05}, it is involved in many conditions associated with decoupling methods in nonlinear analysis and optimization.
The earlier publications \cite{BorweinZhu96,BorweinIoffe1996} employ also a simplified version of \eqref{La}:
\begin{align}
\label{La0}
\Lambda_U^\circ(\varphi_1,\varphi_2):=& \liminf_{\substack{d(x_1,x_2)\to0;\, x_1,x_2\in U }} (\varphi_1(x_1)+\varphi_2(x_2)).
\end{align}
As shown in \cref{P3.1}\,\ref{P3.1.6}, definitions \eqref{La} and \eqref{La0} are not too different,
especially in the situation of our main interest in the current paper when $U$ represents a \nbh\ of a point in many situations.

It follows directly from definitions \eqref{La} and \eqref{La0} that
\[
\Lambda_U(\varphi_1,\varphi_2)\le \Lambda_U^\circ(\varphi_1,\varphi_2)\le \inf_{U}(\varphi_1+\varphi_2),
\]
and the inequalities can be strict (see \cref{E3.3}).
The requirements that $\Lambda_U(\varphi_1,\varphi_2)$ or $\Lambda_U^\circ(\varphi_1,\varphi_2)$
coincide with the conventional infimum of $\varphi_1+\varphi_2$ represent important \emph{qualification conditions}.
If
\begin{gather}
\label{La0qc}
\inf_{U} (\varphi_1+\varphi_2)\le
\Lambda_U^\circ(\varphi_1,\varphi_2)
\end{gather}
(in view of the above, it can only hold as equality),
then the pair $(\varphi_1,\varphi_2)$ is said to be \emph{uniformly lower semicontinuous} \cite[Definition~2.6]{BorweinZhu96} on $U$ (see also
\cite[Remark~2]{BorweinIoffe1996} and \cite[Section~2.3]{Las01}) or \emph{quasicoherent} \cite[Lemma~1.124]{Pen13}.
Some sufficient conditions for this property can be found in \cite{Las01,Pen13}; see \cref{S4}.

A more restrictive sequential definition of uniform lower semicontinuity
(\emph{(ULC) condition}) was introduced in \cite[Definition~6]{BorweinIoffe1996}
(see also \cite[Definition~3.3.17]{BorZhu05}):
$(\varphi_1,\varphi_2)$ is \emph{sequentially uniformly lower semicontinuous} (or \emph{coherent} \cite[Lemma~1.124]{Pen13}) on $U$
if, for any sequences $\{x_{1k}\},\{x_{2k}\}\subset U$
satisfying $d(x_{1k},x_{2k})\to0$ as $k\to+\infty$,
there exists a sequence $\{x_k\}\subset U$ such that
\begin{subequations}\label{ULC}
\begin{align}
\label{ULCa}
&\lim_{k\to+\infty} d(x_k,x_{1k})= \lim_{k\to+\infty} d(x_k,x_{2k})=0,
\\
\label{ULCb}
&\limsup_{k\to+\infty} \big((\varphi_1+\varphi_2)(x_k)- \varphi_1(x_{1k})-\varphi_2(x_{2k})\big)\le0.
\end{align}
\end{subequations}
This definition was formulated in \cite{BorweinIoffe1996,BorZhu05} for the case when $U$ is a ball in a Banach space,
but is meaningful in our more general setting, too.
At the same time, one needs to be a little more careful to ensure that the expression under the $\limsup$ in \eqref{ULCb} is well defined.
It suffices to assume that $\{x_{1k}\}\subset\dom\varphi_1$ and $\{x_{2k}\}\subset\dom\varphi_2$.
The key point that distinguishes this definition from the one in the previous paragraph is the presence of conditions \eqref{ULCa},
which relate the variable of $\varphi_1+\varphi_2$ with those of the decoupled sum $(x_1,x_2)\mapsto\varphi_1(x_1)+\varphi_2(x_2)$.
Recall that the minimizing sequences involved in the expressions compared in \eqref{La0qc} (see \eqref{La0}) are entirely independent.
As observed in \cite[Section~3.3.8]{BorZhu05} (see also \cref{P4.4}), sequential uniform lower semicontinuity possesses certain stability
which makes it more convenient in applications.

Thanks to \cref{P4.3}\,\ref{P4.3.3}
the sequential uniform lower semicontinuity property admits an equivalent analytical representation.
We call it
\emph{firm uniform lower semicontinuity}; see \cref{D4.1}\,\ref{D4.1.3}.

With the problem of minimizing of $\varphi_1+\varphi_2$ in mind,
employing the uniform infimum \eqref{La} naturally leads to the definition of \emph{local uniform minimum} \cite[Section~2.2]{Las01}:
\begin{align}
\label{UM}
(\varphi_1+\varphi_2)(\bx)= \Lambda_{B_\de(\bx)}(\varphi_1,\varphi_2)
\qdtx{for some}\de>0
\end{align}
(see also
\cite[formula~(3.3.2)]{BorZhu05}).
Here $B_\de(\bx)$ stands for the open ball with center at a local minimum $\bx$ and radius $\de>0$.
Expression $\Lambda_{B_\de(\bx)}(\varphi_1,\varphi_2)$ can be replaced in \eqref{UM} with
$\Lambda_{B_\de(\bx)}^\circ(\varphi_1,\varphi_2)$.
If instead of $B_\de(\bx)$ an arbitrary subset $U\subset X$ containing $\bx$ is used, we talk about a {uniform minimum} (robust minimum \cite{Iof12}) of $\varphi_1+\varphi_2$ on $U$.
Every local uniform minimum is obviously a conventional local minimum.
Thus, local minimality conditions established using the decoupling approach
contain either qualification conditions of the type \eqref{UM} explicitly or some sufficient conditions ensuring their fulfillment.
Condition \eqref{UM} is satisfied, in particular, if $\bx$ is a local minimum,
and $(\varphi_1,\varphi_2)$ is sequentially uniformly \lsc\ \cite[Exercise~3.3.5\,(i)]{BorZhu05}.
Some other typical sufficient conditions are collected in \cite[Proposition~2.3]{Las01} and \cite[Proposition~3.3.2]{BorZhu05}.
As emphasized in \cite[Section~3.3.8]{BorZhu05}, without assumption \eqref{UM} some optimality conditions may fail,
while condition \eqref{UM} itself is not tight.

Employing the decoupled definitions \eqref{La} and \eqref{La0}
and the respective associated concepts of (firm) uniform lower semicontinuity and local uniform minimum
allows one to streamline proofs of optimality conditions and calculus relations,
unify and simplify the respective statements,
as well as clarify and in many cases weaken the assumptions.
For instance, it was emphasized in \cite[Remark~2]{BorweinIoffe1996} that (ULC) condition
(firm uniform lower semicontinuity in the language adopted in the current paper)
covers the three types of situations in which (strong) \emph{fuzzy} calculus rules had been established for appropriate subdifferentials in Banach spaces:
when the underlying space is finite-dimensional, when one of the functions has compact level sets
and when all but one functions are uniformly continuous.

Among the fuzzy calculus rules the following (strong) \emph{fuzzy sum rule} is central:
\smallskip

\fbox{\it
\begin{minipage}{0.9\linewidth}
For any $\bx\in\dom\varphi_1\cap\dom\varphi_2$, $x^*\in{\sdf}(\varphi_1+\varphi_2)(\bar x)$ and $\varepsilon>0$,
there exist points $x_1,x_2$
such that
\begin{subequations}
\label{SR}
\begin{gather}
\label{SR-1}
\norm{x_i-\bar x}<\varepsilon,\quad |\varphi_i(x_i)-\varphi_i(\bx)|<\varepsilon\quad
(i=1,2),
\\
\label{SR-2}
\dist(x^*,{\sdf}\varphi_1(x_1)+{\sdf}\varphi_2(x_2)) <\varepsilon.
\end{gather}
\end{subequations}
\end{minipage}
}\smallskip
\\
Here, $\sd$ usually stands for the \Fr\ subdifferential.
This type of rules have been established in appropriate spaces also for other subdifferentials; see
\cite{BorweinZhu96,BorweinIoffe1996,Las01,BorZhu05}.

Note that none of the aforementioned three types of situations involves the traditional (for this type of results) assumption
that all but one functions are locally Lipschitz continuous, thus, ruining the widely spread (even now) myth
that Lipschitzness is absolutely necessary, at least, in infinite-dimensional spaces.

The fact that in finite dimensions the above fuzzy sum rule is valid for arbitrary \lsc\ functions
has been known since the mid-1980s; see \cite[Theorem~2]{Ioffe84}.
A similar result is true also for weakly \lsc\ functions in Hilbert spaces; it is usually formulated in terms of \emph{proximal subdifferentials};
see \cite[Theorem~1.8.3]{ClaLedSteWol98}.
By means of an example, it has been shown in \cite{VanderwerffZhu98}
that the Hilbert space fuzzy sum rule fails if the weak sequential lower semicontinuity is replaced by just lower semicontinuity.
Both the finite-dimensional and Hilbert space fuzzy sum rules can be proved without using the \EVP.
In more general spaces some additional assumptions are required like compactness of the level sets of one of the functions
or uniform (but not necessarily Lipschitz) continuity of all but one functions.
The decoupling approach formalized in \cite{BorweinZhu96,BorweinIoffe1996,Las01,BorZhu05,Pen13}
allows one to treat all these situations within the same framework.
Note that, unlike the finite-dimensional case, in infinite dimensions without additional assumptions strong fuzzy sum rules may fail;
see a counterexample in \cite[Theorem~1]{VanderwerffZhu98}.
For \Fr\ subdifferentials, even with the mentioned additional assumptions
such a rule is only valid in Asplund spaces, and this property is characteristic of Asplund spaces; see \cite{Mor06.1}.

In contrast to the sum rule above, the so-called \emph{weak fuzzy sum rule} is valid  for \lsc\ functions in (appropriate)
infinite-dimensional spaces without additional assumptions; see \cite{Ioffe83,BorweinZhu96,Iof00,BorZhu05,Pen13}.
Instead of condition \eqref{SR-2} involving the distance, it employs the condition
\begin{gather}
\label{WSR}
x^*\in{\sdf} \varphi_1(x_1) +{\sdf} \varphi_2(x_2)+U^*,
\end{gather}
where $U^*$ is an arbitrary weak* neighbourhood of zero in the dual space.
The weak fuzzy sum rule immediately yields the validity of the conventional (strong) fuzzy sum rule in finite dimensions.

The decoupling approach has proved to be useful also when developing sequential subdifferential calculus without constraint qualifications in the convex setting; see Thibault \cite{Thi97}.
Links between the decoupling concepts discussed above and penalization in optimization are discussed in \cite{Ioffe83,Ioffe84,ClaLedSteWol98,Iof12,Pen13}.
Interesting adaptations of the notions of uniform infimum and uniform minimum for characterizing the subdifferential of the supremum of an infinite family of functions have been suggested recently by P\'{e}rez-Aros \cite{Per19,Per19.2}.

In our recent paper \cite{KruMeh22} the decoupling approach was used intuitively when proving the main result \cite[Theorem~4.1]{KruMeh22}.
When analyzing later the proof of that theorem and related definitions and facts, and tracing the ideas back to the foundations in
\cite{BorweinZhu96,BorweinIoffe1996,Las01,BorZhu05}, we have realized that the `novel notions of semicontinuity'
discussed in \cite[Section~3]{KruMeh22} are closely related to the uniform lower semicontinuity as in \cite{BorweinIoffe1996,BorweinZhu96,Las01}.
More importantly, our version of uniform lower semicontinuity is actually weaker,
thus, leading to weaker notions of uniform infimum, firm uniform infimum and local uniform minimum
as well as fuzzy optimality conditions and calculus relations under weaker assumptions.
Further developing the notions introduced and studied in \cite{KruMeh22} and putting them into the context of the general theory
developed in \cite{BorweinIoffe1996,BorweinZhu96,Las01,BorZhu05,Pen13} is the main aim of the current paper.
Along the way, we unify the terminology and notation, and fill in some gaps in the general theory.

We clearly distinguish between the uniform lower semicontinuity defined by \eqref{La0qc} and the firm uniform lower semicontinuity
(the analytical counterpart of the sequential lower semicontinuity defined using \eqref{ULC}; see \cref{D4.1}\,\ref{D4.1.3})
exposing the advantages of the latter stronger property.
The novel weaker properties arising from \cite{KruMeh22} are called \emph{quasiuniform lower semicontinuity} and \emph{firm quasiuniform lower semicontinuity}.
The first one is defined similarly to \eqref{La0qc} using instead of \eqref{La0}  the \emph{quasiuniform infimum}
\begin{gather}
\label{M}
{\Lambda}_U^\dag(\varphi_1,\varphi_2)
:=
\inf_{V\in EI(U)}\; \Lambda_V^\circ(\varphi_1,\varphi_2)
\end{gather}
of $(\varphi_1,\varphi_2)$ over $U$.
Here $EI(U)$ stands for the collection of all \emph{essentially interior} subsets of $U$; see \cref{D2.1}.
Clearly, ${\Lambda}_U^\dag(\varphi_1,\varphi_2)\ge {\Lambda}_U^\circ(\varphi_1,\varphi_2)$, and the inequality can be strict; see \cref{E3.2,E3.3}.
To simplify the comparison, all four uniform lower semicontinuity notions together with their localized (near a point) versions are collected in a single \cref{D4.1}.

We study the weaker than \eqref{UM} local \emph{quasiuniform minimality} notion:
\begin{align}
\label{QUM}
(\varphi_1+\varphi_2)(\bx)= \Lambda_{B_\de(\bx)}^\dagger(\varphi_1,\varphi_2)
\qdtx{for some}\de>0,
\end{align}
employing the quasiuniform infimum \eqref{M},
together with the related notions of \emph{quasiuniform stationarity} and \emph{quasiuniform $\eps$-minimality}; see \cref{D6.1}.
Using these new notions allows one to formulate more subtle conditions.
The mentioned quasiuniform minimality/stationarity coincide with the corresponding conventional local minimality,
stationarity and $\eps$-minimality when the pair $(\varphi_1,\varphi_2)$ is quasiuniformly lower semicontinuous on an appropriate \nbh\ of $\bx$.
We establish rather general primal and dual (fuzzy multiplier rules) necessary conditions characterizing quasiuniform
$\eps$-minimum of the sum of two functions.
Under the assumption of quasiuniform lower semicontinuity of $(\varphi_1,\varphi_2)$,
they characterize the conventional $\eps$-minimum and, as a consequence, also any stationary point and local minimum.
The sufficient conditions for quasiuniform lower semicontinuity discussed in the paper encompass all known conditions ensuring fuzzy multiplier rules.
These are formulated in general Banach spaces in terms of Clarke subdifferentials and in Asplund spaces in terms of \Fr\ subdifferentials.
In general Banach spaces Clarke subdifferentials can be replaced in this type of statements by the \emph{$G$-subdifferentials} of Ioffe \cite{Iof17}.
The mentioned fuzzy multiplier rules naturally lead to certain fuzzy subdifferential calculus results under rather weak assumptions; see \cref{T7.1}.

The structure of the paper is as follows.
The next \cref{S2} contains the basic notation and some preliminary results from variational analysis used throughout the paper.
In particular, we introduce essentially interior subsets that are key for the definition of quasiuniform infimum
and other new notions, and can be useful elsewhere.
In \cref{S3} we discuss the notions of uniform and quasiuniform infimum
and two other `decoupling quantities' as well as several analogues of the qualification condition \eqref{La0qc}.
Diverse examples illustrate the computation of the `decoupling quantities'.
These `decoupling quantities' and qualification conditions provide the basis for the definitions of uniform and quasiuniform lower semicontinuity
and their `firm' analogues discussed in \cref{S4}.
We show that firm uniform and firm quasiuniform lower semicontinuity properties are stable under
uniformly continuous perturbations of the involved functions and prove several sufficient conditions
for the mentioned uniform lower semicontinuity properties.
In \cref{S5} we investigate the situation where at least one of the involved functions is the indicator function
of a set and discuss the notions of relative uniform and quasiuniform lower semicontinuity.
We show that the situations when a
pair of indicator functions are not firmly uniformly or firmly quasiuniformly lower semicontinuous are rare.
The firm quasiuniform lower semicontinuity of a
pair of indicator functions near a point in the intersection of the sets is implied, for instance,
by the well-known and widely used subtransversality property.
\cref{S6} is devoted to the problem of minimizing the sum of two functions.
Here we prove rather general primal and dual necessary conditions characterizing quasiuniform
$\eps$-minimum and formulate several consequences.
In \cref{S7} we illustrate the value of quasiuniform lower semicontinuity in the context
of subdifferential calculus.
An application in sparse optimal control is considered in \cref{S8}.
The paper closes with some concluding remarks in \cref{sec:conclusions}.

\section{Notation and preliminaries}
\label{S2}

\paragraph{Basic notation and definitions}

Our basic notation is standard; cf. e.g.
\cite{Iof17,KlaKum02,Mor06.1}.
Throughout the paper $X$ and $Y$ are either metric or normed spaces (or more specifically Banach or Asplund spaces).
We use the same notation
$d(\cdot,\cdot)$ and $\norm{\cdot}$
for distances and norms in all spaces
(possibly with a subscript specifying the space).
Normed spaces are often treated as metric spaces with the distance determined by the norm.
If $X$ is a normed space, $x\in X$ and $U\subset X$,
we use $x+U:=U+x:=\{x+u\,|\,u\in U\}$ for brevity of notation.
The topological dual of $X$ is denoted by $X^*$, while $\langle\cdot,\cdot\rangle$
denotes the bilinear form defining the pairing between the two spaces.
If not stated otherwise, product spaces are equipped with the associated
maximum distances or norms.
The associated dual norm is the sum norm.
In a metric space, $\B$ and $\overline\B$ are the \emph{open} and \emph{closed} unit balls, while $B_\de(x)$ and $\overline B_\de(x)$ are the
\emph{open} and \emph{closed}
balls with radius $\de>0$ and center $x$, respectively.
We write $\B^*$ and $\overline\B^*$ to denote the \emph{open} and \emph{closed} unit balls in the dual to a normed space.

The distance from a point $x$ to a set $U$ is defined
by $\dist(x,U):=\inf_{x\in U}d(x,u)$ with the convention $\dist(x,\es):=+\infty$.
Furthermore, for two sets $U_1,U_2$, we define $\dist(U_1,U_2):= \inf_{x\in U_1}\dist(x,U_2)$.
In this paper, $\Int U$ and $\cl U$ represent the interior and the closure of
$U$, respectively.
In a normed space we use $\co U$ and $\cl\co U$ to denote
the convex hull and the closed convex hull of $U$, respectively.
We write $x_k\to\bar x$ to denote the (strong) convergence of a sequence $\{x_k\}$ to a point $\bar x$.
In a normed space, $x_k\weakly\bar x$ expresses the weak convergence of
$\{x_k\}$ to $\bar x$, i.e., $\langle x^*,x_k\rangle\to\langle x^*,\bar x\rangle$
for each $x^*\in X^*$.
Similarly, $x_k^*\overset{*}{\weakly}x^*$ denotes the weak* convergence of
$\{x_k^*\}$ to $x^*$ in the dual space.
\sloppy

Symbols $\R$, $\R_+$ and $\N$ represent the sets of all real numbers, all nonnegative real numbers and all positive integers, respectively.
We make use of the notation $\R_\infty:=\R\cup\{+\infty\}$ and the conventions
$\inf\es_{\R}=+\infty$ and
$\sup\es_{\R_+}=0$, where $\es$ (possibly with a subscript) denotes the empty subset (of a given set).

\begin{definition}
\label{D2.1}
Let $X$ be a metric space and $U\subset X$.
A set $V\subset X$ is an \emph{essentially interior} subset of $U$ if $B_\rho(V):=\bigcup_{x\in V}B_\rho(x)\subset U$ for some $\rho>0$, or equivalently, if
$\inf_{x\in V}\dist(x,X\setminus U)>0$.
We write $EI(U)$ and $\EI(U)$ to denote the collections of, respectively, all essentially interior subsets and all {closed} essentially interior subsets of $U$.
\end{definition}

The next lemma summarizes basic properties of essentially interior subsets of a given set.

\begin{lemma}
\label{L2.1}
Let $X$ be a metric space and $U\subset X$.
The following assertions hold:
\begin{enumerate}
\item
$\es\in\EI(U)$ and $X\in\EI(X)$.
\item
\label{L2.1.2}
Let $V\in EI(U)$.
If $V'\subset V$ and $U\subset U'$, then $V'\in EI(U')$.
In particular, $\Int V\in EI(U)$ and $V\in EI(\cl U)$.
\item
\label{L2.1.3}
If $x\in\Int U$, then $\{x\}\in\EI(U)$.
\item
\label{L2.1.4}
$\bigcup_{V\in EI(U)}V=\bigcup_{V\in\EI(U)}V=\Int U$.
\item
\label{L2.1.5}
If $V\in EI(U)$, then
$\cl V\in\EI(U)$.
\item
\label{L2.1.6}
If $V\in EI(U)$, then there exists a subset $V'\in EI(U)$ such that $V\in EI(V')$.
\item
\label{L2.1.7}
$EI(U)=EI(\Int U)$ and
$\EI(U)=\EI(\Int U)$.
\item
\label{L2.1.8}
Let $\bx\in X$ and $\de>0$.
Then $B_{\rho}(\bx)\in EI(B_\de(\bx))$ for all $\rho\in(0,\de)$.
\cnta
\end{enumerate}
Suppose $X$ is a normed space.
\begin{enumerate}
\cntb
\item
\label{L2.1.9}
Let $\bx\in X$ and $\de>0$.
If $V\in EI(\overline B_\de(\bx))$, then $V\subset B_{\rho}(\bx)$ for some $\rho\in(0,\de)$.
\item
\label{L2.1.10}
If $U$ is convex and $V\in EI(U)$, then
$\co V\in EI(U)$.
\end{enumerate}
\end{lemma}

\begin{proof}
Most of the assertions are direct consequences of \cref{D2.1}.
\begin{enumerate}
\item
Both $\es$ and $X$ are closed, and we have $B_\rho(\es)=\es\subset U$ and $B_\rho(X)=X$ for any $\rho>0$.
\item
Since $V\in EI(U)$, we have $B_\rho(V)\subset U$ for some $\rho>0$.
If $V'\subset V$ and $U\subset U'$, then
$B_\rho(V')\subset B_\rho(V)\subset U\subset U'$; hence $V'\in EI(U')$.
\item
If $x\in\Int U$, then $B_\rho(\{x\})= B_\rho(x)\subset U$ for some $\rho>0$, and $\{x\}$ is closed.
\item
Observe that $\EI(U)\subset EI(U)$; hence $\bigcup_{V\in\EI(U)}V\subset\bigcup_{V\in EI(U)}V$.
Next, if $x\in V\in EI(U)$, then $B_\rho(x)\subset U$ for some $\rho>0$; hence,
$x\in\Int U$.
Finally, if $x\in\Int U$, then, by \ref{L2.1.3}, $x\in\{x\}\in\EI(U)$.

\item
If $V\in EI(U)$, then,
thanks to the continuity of the distance function, we have
$$\inf_{x\in\cl V}\dist(x,X\setminus U)=\inf_{x\in V}\dist(x,X\setminus U)>0.$$
\item
Let $V\in EI(U)$, i.e., $B_\rho(V)\subset U$ for some $\rho>0$.
Then the set $V':=B_{\rho/2}(V)$ validates the assertion.
\item
By \ref{L2.1.2}, $EI(\Int U)\subset EI(U)$.
Conversely, let $V\in EI(U)$, i.e., $B_\rho(V)\subset U$ for some $\rho>0$.
Then $B_{\rho/2}(V)\subset\Int U$, i.e., $V\in EI(\Int U)$.
The above argument applies, in particular, to closed essentially interior subsets.
\item
If $\rho\in(0,\de)$, then, for any $\rho'\in(0,\de-\rho)$, we have $B_{\rho'}(B_\rho(\bx))\subset B_\de(\bx)$.
\item
Let $V\in EI(\overline B_\de(\bx))$, i.e.,
$B_{r}(V)\subset \overline B_\de(\bx)$
for some $r>0$.
Then $\overline B_r(V)\subset\overline B_\de(\bx)$.
Set $\rho:=\de-r/2<\de$.
Let $x\in V$ and $x\ne\bx$.
We are going to show that $\|x-\bx\|<\rho$.
Set $\tilde x:=x+r(x-\bar x)/\norm{x-\bar x}$.
Then $\tilde x\in\overline B_r(V)\subset\overline B_\de(\bx)$, and consequently, $\norm{\tilde x-\bar x}\le\de$.
Observe that $\norm{\tilde x-\bar x}=\norm{x-\bar x}+r$.
Hence, $\|x-\bx\|\le\de-r<\de-r/2=\rho$, and consequently, $\rho>0$ and $V\subset B_{\rho}(\bx)$.
\item
Let $V\in EI(U)$.
By \cref{D2.1}, $V+\rho\B\subset U$ for some $\rho>0$.
Let $x\in\co V$ and $x'\in\rho\B$.
Then $x=\sum_{i=1}^n\lambda_i x^i$ for some $n\in\N$, $x^1,\ldots,x^n\in V$ and $\lambda_1,\ldots,\lambda_n\in\R_+$ with $\sum_{i=1}^n\lambda_i=1$, and consequently, $x+x'=\sum_{i=1}^n\lambda_i x^i+x'=\sum_{i=1}^n\lambda_i(x^i+x')\in U.$
Hence, $\co V+\rho\B\subset U$, i.e.,
$\co V\in EI(U)$.
\end{enumerate}
\end{proof}

\begin{remark}
Assertion \ref{L2.1.9} of \cref{L2.1}
may fail if $X$ is merely a metric space.
Indeed,
let $X$ be the closed interval $[0,2]$ in $\R$ with the conventional distance,
$\bar x:=1$ and $V:=\{0\}$.
It is easy to see that $B_{2}(V)\subset \overline B_1(\bx)=X$; hence $V\in EI(\overline B_1(\bx))$.
However, $V\not\subset B_{\rho}(\bx)$ for any $\rho\in(0,1)$.
Note that in this example
$B_1(\bx)\ne\Int\overline B_1(\bx)=\overline B_1(\bx)$ and
$V\notin EI(B_1(\bx))$.
\end{remark}

For an extended-real-valued function $f\colon X\to\R_\infty$,
its domain
and epigraph are defined,
respectively,
by
$\dom f:=\{x \in X\mid {f(x) < +\infty}\}$
and
$\epi f:=\{(x,\alpha) \in X \times \mathbb{R}\mid {f(x) \le \alpha}\}$.
If $X$ is a metric space,
a point $\bar x\in\dom f$ is called \emph{stationary} for $f$
whenever
\[
\liminf\limits_{x\to\bar x,\,x\neq\bar x}\frac{f(x)-f(\bar x)}{d(x,\bar x)}\geq0.
\]
Clearly, if $\bar x$ is a local minimum of $f$, then it is stationary.
If $\bar x\in U\subset X$ and
$\inf_{U}f>f(\bx)-\eps$ for some $\varepsilon>0$, then $\bx$ is an \emph{$\eps$-mini\-mum} of $f$ on $U$.
The next lemma provides a connection between stationarity and $\varepsilon$-minimality.

\begin{lemma}
\label{lem:characterization_stationarity}
	Let $X$ be a metric space, $f\colon X\to\R_\infty$ and $\bar x\in\dom f$.
	Then $\bar x$ is stationary for $f$ if and only if
	for any $\varepsilon>0$, there is a $\de_\varepsilon>0$ such that, for any $\de\in(0,\de_\varepsilon)$,
	$\bar x$ is an $\varepsilon\de$-minimum of $f$ on $B_\de(\bar x)$.
\end{lemma}

\begin{proof}
	If $\bar x$ is stationary for $f$,
	then for any $\varepsilon>0$, there is a $\de_\varepsilon>0$ such that
	$f(x)-f(\bar x)> -\varepsilon d(x,\bar x)$ for all $x\in\overline B_{\de_\varepsilon}(\bar x)$.
	Particularly, for all $\de\in(0,\de_\varepsilon)$ and $x\in \overline B_{\de}(\bar x)$, we have
	$f(x)-f(\bar x)> -\varepsilon\de$.
	Conversely, suppose that $\bar x$ is not stationary for $f$.
	Then there is an $\eps>0$ such that
	$\liminf_{x\to\bar x,\,x\neq\bar x}(f(x)-f(\bar x))/d(x,\bar x)<-\eps$.
	Choose any number $\de_\varepsilon>0$.
	Then there exists a point $\hat x$ such that $d(\hat x,\bx)<\de_\varepsilon$ and $f(\hat x)-f(\bar x)<-\varepsilon d(\hat x,\bar x)$,
	and we can find a number $\rho>0$ such that $d(\hat x,\bx)+\rho<\de_\varepsilon$ and $f(\hat x)-f(\bar x)<-\varepsilon(d(\hat x,\bar x)+\rho)$.
	Set $\de:=d(\hat x,\bar x)+\rho$ and observe that $\hat x\in B_\de(\bar x)$ and $f(\hat x)-f(\bar x)<-\varepsilon\de$,
	i.e., $\bar x$ is not an $\varepsilon\de$-minimum of $f$ on $B_\de(\bar x)$.
\end{proof}

For a set-valued mapping $F\colon X\toto Y$,
its domain
and graph are defined
respectively,
by
$\dom F:=\{x\in X\mid F(x)\ne\es\}$
and
$\gph F:=\{(x,y) \in X \times Y\mid y\in F(x)\}$.

Recall that a Banach space is \emph{Asplund} if every continuous convex function on an open convex set is Fr\'echet differentiable on a dense subset \cite{Phe93}, or equivalently, if the dual of each
separable subspace is separable.
A Banach space is \emph{\Fr\ smooth} if it has an equivalent norm that is \Fr\ differentiable away from zero \cite{KruMor80,Kru03,BorZhu05}.
All reflexive, particularly, all finite-dimensional Banach spaces are \Fr\ smooth, while all \Fr\ smooth spaces are Asplund.
We refer the reader to
\cite{DevGodZiz93.2,Phe93,BorZhu05,Mor06.1} for discussions about and characterizations of Asplund and \Fr\ smooth spaces.

\paragraph{Subdifferentials, normal cones and coderivatives}

Below we review some standard notions of generalized differentiation
which can be found in many monographs; see, e.g., \cite{Cla83,Iof17,Mor06.1}.
Below $X$ and $Y$ are normed spaces.

For a function $\varphi\colon X\to\R_\infty$ and a point $\bar x\in\dom\varphi$, the (possibly empty) set
\[
	\sd \varphi(\bar x)
	:=
	\left\{x^*\in X^*\,\middle|\,
		\liminf_{x\to\bar x,\,x\neq\bar x}
			\frac{\varphi(x)-\varphi(\bar x)-\langle x^*,x-\bar x\rangle}{\norm{x-\bar x}}
		\geq 0
	\right\}
\]
is the \emph{Fr\'{e}chet subdifferential} of $\varphi$ at $\bar x$.
If $\bx$ is a local minimum (or, more generally, a stationary point) of $\varphi$, then obviously $0\in\sd\varphi(\bx)$ (\emph{Fermat rule}).
If $X$ is Asplund, the \emph{limiting subdifferential} of $\varphi$ at $\bar x$ can be defined as
the outer/upper limit (with respect to the norm topology in $X$ and weak* topology in $X^*$) of Fr\'{e}chet subdifferentials:
\begin{align*}
	\sdm\varphi(\bar x)
	:=
	\left\{
		x^*\in X^*\,\middle|\,
		\begin{aligned}
			&\exists\{x_k\}\subset X,\,\{x_k^*\}\subset X^*\colon\\
			&\qquad
			x_k\to\bar x,\,\varphi(x_k)\to\varphi(\bar x),\,
			x_k^*\overset{*}{\weakly}x^*,\, x_k^*\in\sd\varphi(x_k)\,\forall k\in\N
		\end{aligned}
	\right\}.
\end{align*}

For a subset $\Omega\subset X$ and a point $\bar x\in\Omega$, the closed convex (possibly trivial) cone
\[
	N_\Omega(\bar x)
	:=
	\left\{
		x^*\in X^*\,\middle|\,
		\limsup_{x\to\bar x,\,x\in\Omega,\,x\neq\bar x}
			\frac{\langle x^*,x-\bar x\rangle}{\norm{x-\bar x}}
		\leq 0
	\right\}
\]
is the \emph{Fr\'{e}chet normal cone} to $\Omega$ at $\bar x$.
It is easy to check that $N_\Omega(\bar x)=\sd i_\Omega(\bar x)$,
where $i_\Omega\colon X\to\R_\infty$ is the \emph{indicator function} of $\Omega$,
given by $i_\Omega(x)=0$ if $x\in\Omega$ and $i_\Omega(x)=+\infty$ otherwise,
and if $\bx\in\dom\varphi$, then
\[
	\sd \varphi(\bar x)
	=
	\left\{x^*\in X^*\,\middle|\,(x^*,-1)\in N_{\epi\varphi}(\bar x,\varphi(\bar x))\right\}.
\]
The \emph{Clarke normal cone} to $\Omega$ at $\bar x$ is defined by
\[
	N^\textup{C}_\Omega(\bar x)
	:=
	\bigl\{x^*\in X^*\,\bigl|\,
		\forall d\in T^\textup{C}_\Omega(\bar x)\colon\,\langle x^*,d\rangle\leq 0
	\bigr\},
\]
where
\[
	T^\textup{C}_\Omega(\bar x)
	:=
	\left\{
		d\in X\,\middle|\,
		\begin{aligned}
		&\forall \{x_k\}\subset \Omega,\,\{t_k\}\subset(0,+\infty)\colon
			\,x_k\to\bar x,\,t_k\downarrow 0\\
		&\quad\Rightarrow\,\exists\{d_k\}\subset X\colon\,d_k\to d,\,x_k+t_kd_k\in\Omega\,
		\forall k\in\N
		\end{aligned}
	\right\}
\]
is the \emph{Clarke tangent cone} to $\Omega$ at $\bar x$, while
the \emph{Clarke subdifferential} of a function $\varphi\colon X\to\R_\infty$ at $\bar x\in\dom\varphi$
can be defined via
\[
	\sdc \varphi(\bar x)
	:=
	\bigl\{x^*\in X^*\,\bigl|\,
		(x^*,-1)\in N^\textup{C}_{\epi\varphi}(\bar x,\varphi(\bar x))
	\bigr\}
\]
(the direct definition of the Clarke subdifferential is a little more involved and employs Clarke--Rockafellar directional derivatives).
It always holds $N_\Omega(\bar x)\subset N^\textup{C}_\Omega(\bar x)$ and $\sd\varphi(\bar x)\subset\sdm\varphi(\bar x)\subset\sdc\varphi(\bar x)$,
and whenever $\Omega$ and $\varphi$ are convex, the above normal cones and subdifferentials reduce to the conventional constructions of convex analysis:
\begin{gather*}
	N_\Omega(\bar x)
	=
	N^\textup{C}_\Omega(\bar x)
	=
	\{x^*\in X^*\,|\,\forall x\in\Omega\colon\,\langle x^*,x-\bar x\rangle\leq 0\},
\\
	\sd\varphi(\bar x)
	=
	\sdm\varphi(\bar x)
	=
	\sdc\varphi(\bar x)
	=
	\{x^*\in X^*\,|\,
		\forall x\in\dom\varphi\colon\,\varphi(x)\geq\varphi(\bar x)+\langle x^*,x-\bar x\rangle
	\}.
\end{gather*}

For a mapping $F\colon X\toto Y$ between normed spaces and a point $(\bx,\by)\in\gph F$,
the set-valued mapping $D^* F(\bx,\by)\colon Y^*\rightrightarrows X^*$ given via
\[
	\forall y^*\in Y^*\colon\quad
	D^*F(\bx,\by)(y^*)
	:=
	\left\{x^*\in X^*\,\middle|\,
		(x^*,-y^*)\in N_{\gph F}(\bar x,\bar y)
	\right\}
\]
is the \emph{Fr\'{e}chet coderivative} of $F$ at $(\bx,\by)$.
If $F$ is single-valued and $\by=F(\bx)$, we write simply $D^* F(\bx)$ for brevity.

\paragraph{Preliminary results}

Here we recall some fundamental results from the literature used in the sequel.
We start with the celebrated Ekeland variational principle; see e.g.\
\cite{AubFra90,KlaKum02,BorZhu05,Mor06.1,Iof17}.

\begin{lemma}
\label{lem:Ekeland}
Let $X$ be a complete metric space,
$\varphi\colon X\to\R_{\infty}$ lower semicontinuous and bounded from below and
$\bx\in\dom\varphi$.
Then, for any $\varepsilon>0$, there exists a point $\hat x\in X$ satisfying the following conditions:
\begin{enumerate}
\item
$\varphi(\hat x)+\varepsilon d(\hat x,\bx)\le \varphi(\bx)$
(as a consequence,
$\varphi(\hat x)\le \varphi(\bx)$);

\item
$\forall x\in X\setminus\{\hat x\}\colon\, \varphi(\hat x)<\varphi(x)+\varepsilon d(x,\hat x)$.
\end{enumerate}
\end{lemma}

The next lemma summarizes some standard sum rules for Fr\'{e}chet and Clarke
subdifferentials as well as conventional subdifferentials of convex functions which can be found
in many monographs on variational
analysis
\cite{Rockafellar1970,IofTik79,Cla83,Phe93,Zal02, BorZhu05,Iof17}.

\begin{lemma}\label{lem:SR}
Let $X$ be a Banach space,
$\varphi_1,\varphi_2\colon X\to\R_\infty$
and $\bar x\in\dom \varphi_1\cap\dom \varphi_2$.
Then the following assertions hold.
\begin{enumerate}
\item
\label{lem:SR.0}
\textbf{Differentiable sum rule} \cite{Kru03}.
If $\varphi_1$ is differentiable at $\bar x$ with derivative $\varphi_1'(\bar x)$, then
${\sdf}(\varphi_1+\varphi_2)(\bar x)= \varphi_1'(\bar x)+{\sdf}\varphi_2(\bar x).$
\item
\label{lem:SR.1}
\textbf{Convex sum rule} \cite{IofTik79,Phe93,Zal02}.
If $\varphi_1$ and $\varphi_2$ are convex, and $\varphi_1$ is continuous at a point in $\dom \varphi_2$, then
${\sdf}(\varphi_1+\varphi_2)(\bar x)={\sdf} \varphi_1(\bar x)+{\sdf} \varphi_2(\bar x).$
\item
\label{lem:SR.Clarke}
\textbf{Clarke sum rule} \cite{Cla83,Roc79}.
If $\varphi_1$ is Lipschitz continuous near $\bar x$,
and $\varphi_2$ is lower semicontinuous near $\bar x$,
then
$\sdc(\varphi_1+\varphi_2)(\bar x)\subset \sdc\varphi_1(\bar x)+\sdc \varphi_2(\bar x).$
\item
\label{lem:SR.2}
\textbf{Fuzzy sum rule} \cite{Fab89,Iof00}.
If $X$ is Asplund, $\varphi_1$ is Lipschitz continuous near $\bar x$,
and $\varphi_2$ is lower semicontinuous near $\bar x$,
then, for any $x^*\in{\sdf}(\varphi_1+\varphi_2)(\bar x)$ and $\varepsilon>0$,
there exist points $x_1,x_2\in X$
such that conditions \eqref{SR} hold true.
\item\label{lem:SR.3}
\textbf{Weak fuzzy sum rule} \cite{Iof00,BorZhu05,Pen13}.
If $X$ is Fr\'{e}chet smooth, and $\varphi_1$ and $\varphi_2$
are lower semicontinuous near $\bar x$, then, for any $x^*\in{\sdf}(\varphi_1+\varphi_2)(\bar x)$, $\varepsilon>0$ and a weak* neighbourhood $U^*$ of zero in $X^*$,
there exist $x_1,x_2\in X$
such that conditions \eqref{SR-1} and \eqref{WSR} hold true.
\end{enumerate}
\end{lemma}

\begin{remark}
\label{R2.1}
\begin{enumerate}
\item
The sum rules in parts
\ref{lem:SR.Clarke} and \ref{lem:SR.2} of \cref{lem:SR}
contain the standard (and commonly believed to be absolutely necessary) assumption
of Lipschitz continuity of one of the functions.
In fact, this assumption is not necessary.
For the fuzzy sum rule in part \ref{lem:SR.2},
it has been shown in \cite[Corollary~3.4\,(ii)]{CutFab16} that it suffices to assume
one of the functions to be \emph{uniformly continuous} in a neighbourhood of the reference point.
In the setting of smooth spaces the latter fact was discussed also in \cite{BorweinZhu96,Las01,BorZhu05,Pen13}.

\item
Part \ref{lem:SR.3} of \cref{lem:SR} shows, in particular, that the fuzzy sum rule holds inherently in finite dimensions without assuming one of the functions to be Lipschitz continuous near the reference point, thus, strengthening the assertion in part \ref{lem:SR.2}.

\item
\label{R2.1.3}
The sum rules in parts \ref{lem:SR.0}, \ref{lem:SR.1} and \ref{lem:SR.Clarke} of \cref{lem:SR} are \emph{exact}
in the sense that the subdifferentials (and the derivative in part \ref{lem:SR.0}) in their \RHS s are computed at the reference point.
In contrast, the rules for Fr\'echet subdifferentials in parts \ref{lem:SR.2} and \ref{lem:SR.3}
are often referred to as \emph{fuzzy} or \emph{approximate} because the subdifferentials in the \RHS s of the inclusions
are computed at some other points arbitrarily close to the reference point.
\end{enumerate}
\end{remark}

\section{Uniform and quasiuniform infimum}
\label{S3}

In this section we discuss the notions of uniform and quasiuniform infimum and two other `decoupling quantities'
as well as several analogues of the qualification condition~\eqref{La0qc}.

Let functions $\varphi_1,\varphi_2\colon X\to\R_\infty$ on a metric space $X$ and a subset $U\subset X$
satisfy
\begin{gather}
\label{nonempty}
\dom\varphi_1\cap\dom\varphi_2\cap U\ne\es.
\end{gather}
Thus, $\inf_{U} (\varphi_1+\varphi_2)<+\infty$.
Recall that the \emph{uniform infimum} of $(\varphi_1,\varphi_2)$ over (or around) $U$ is defined by either \eqref{La} or \eqref{La0}, while the \emph{quasiuniform infimum} ${\Lambda}_U^\dag(\varphi_1,\varphi_2)$ of $(\varphi_1,\varphi_2)$ over $U$ is defined by \eqref{M}.
Some elementary properties of these quantities are collected in the next proposition.

\begin{proposition}
\label{P3.1}
\begin{enumerate}
\item
\label{P3.1.1}
$\Lambda_U(\varphi_1,\varphi_2)\le \Lambda_U^\circ(\varphi_1,\varphi_2)\le
\inf_{U} (\varphi_1+\varphi_2)$.
\item
\label{P3.1.2}
$\Lambda_X(\varphi_1,\varphi_2)= \Lambda_X^\circ(\varphi_1,\varphi_2)= {\Lambda}_X^\dag(\varphi_1,\varphi_2)$.
\item
\label{P3.1.3}
If $U_1\subset U_2$, then
$\Lambda_{U_1}(\varphi_1,\varphi_2)\ge
\Lambda_{U_2}(\varphi_1,\varphi_2)$,
$\Lambda_{U_1}^\circ(\varphi_1,\varphi_2)\ge
\Lambda_{U_2}^\circ(\varphi_1,\varphi_2)$
and
$\Lambda_{U_1}^\dag(\varphi_1,\varphi_2)\ge
\Lambda_{U_2}^\dag(\varphi_1,\varphi_2)$.
\sloppy
\item
\label{P3.1.4}
$\Lambda_U^\circ(\varphi_1,\varphi_2)\le {\Lambda}_U^\dag(\varphi_1,\varphi_2)\le
\inf_{\Int U} (\varphi_1+\varphi_2)$.
\item
\label{P3.1.5}
If $\Int U=\es$, then ${\Lambda}_U^\dag(\varphi_1,\varphi_2)=+\infty$.
\item
\label{P3.1.6}
$\Lambda_U(\varphi_1,\varphi_2)= \lim_{\eps\downarrow0} \Lambda_{B_\eps(U)}^\circ(\varphi_1,\varphi_2)$.
\item
\label{P3.1.7}
The following representations hold true:
\begin{align*}
{\Lambda}_U^\dag(\varphi_1,\varphi_2)&= \inf\limits_{V\in EI(U)}\; \liminf\limits_{\substack{d(x_1,x_2)\to0,\; x_1\in V}} (\varphi_1(x_1)+\varphi_2(x_2))
\\&=
\inf\limits_{V\in\EI(U)}\; \liminf\limits_{\substack{d(x_1,x_2)\to0,\; x_1\in V}} (\varphi_1(x_1)+\varphi_2(x_2))
=\inf\limits_{V\in\EI(U)}\; 
{\Lambda}_V^\circ(\varphi_1,\varphi_2). \end{align*}

\item
\label{P3.1.9}
If $X$ is a normed space, $\bx\in X$ and $\de>0$, then the following representations hold true:
\begin{align*}
{\Lambda}_{B_\de(\bx)}^\dag(\varphi_1,\varphi_2)= {\Lambda}_{\overline B_\de(\bx)}^\dag(\varphi_1,\varphi_2)&= \inf_{\rho\in(0,\de)}\; {\Lambda}_{B_\rho(\bx)}^\circ(\varphi_1,\varphi_2)= \inf_{\rho\in(0,\de)}\; {\Lambda}_{\overline B_\rho(\bx)}^\circ(\varphi_1,\varphi_2)\\
&= \inf_{\rho\in(0,\de)}\; \liminf_{\substack{d(x_1,x_2)\to0,\;
x_1\in B_\rho(\bx)}} (\varphi_1(x_1)+\varphi_2(x_2))\\
&= \inf_{\rho\in(0,\de)}\; \liminf_{\substack{d(x_1,x_2)\to0,\;
x_1\in\overline B_\rho(\bx)}} (\varphi_1(x_1)+\varphi_2(x_2)).
\end{align*}
\end{enumerate}
\end{proposition}

\begin{proof}
Assertions \ref{P3.1.1}, \ref{P3.1.2} and \ref{P3.1.3} are immediate consequences of definitions \eqref{La}, \eqref{La0} and \eqref{M}.
Assertions \ref{P3.1.4} and \ref{P3.1.5} also take into account \cref{L2.1}\,\ref{L2.1.4} and the convention
$\inf\es_{\R}=+\infty$, while assertion \ref{P3.1.9} is a consequence of definition \eqref{M}, assertion \ref{P3.1.7},
and parts \ref{L2.1.2}, \ref{L2.1.5}, \ref{L2.1.8} and \ref{L2.1.9} of \cref{L2.1}.
We prove assertions \ref{P3.1.6} and \ref{P3.1.7}.

\ref{P3.1.6}.
It follows from \eqref{La} and \eqref{La0} that $\Lambda_{B_\eps(U)}^\circ(\varphi_1,\varphi_2)\le \Lambda_U(\varphi_1,\varphi_2)$ for all $\eps>0$.
In view of \ref{P3.1.3}, $\Lambda_{B_\eps(U)}^\circ(\varphi_1,\varphi_2)$ is a nonincreasing function of $\eps$.
Hence,
$\lim_{\eps\downarrow0} \Lambda_{B_\eps(U)}^\circ(\varphi_1,\varphi_2)\le \Lambda_U(\varphi_1,\varphi_2)$.
Let $\lim_{\eps\downarrow0} \Lambda_{B_\eps(U)}^\circ(\varphi_1,\varphi_2)<\ga$ and $\eta>0$.
Then there is an $\eps\in(0,\eta)$ such that $\Lambda_{B_\eps(U)}^\circ(\varphi_1,\varphi_2)<\ga$.
By definition \eqref{La0}, there exist $x_1,x_2\in B_\eps(U)$ such that $d(x_1,x_2)<\eta$ and $\varphi_1(x_1)+\varphi_2(x_2)<\ga$.
Thus, $\dist(x_1,U)<\eps<\eta$.
As $\eta>0$ is an arbitrary number, it follows from definition \eqref{La} that $\Lambda_U(\varphi_1,\varphi_2)\le\ga$,
and consequently, $\Lambda_U(\varphi_1,\varphi_2)\le \lim_{\eps\downarrow0} \Lambda_{B_\eps(U)}^\circ(\varphi_1,\varphi_2)$.
Combining both estimates gives~\ref{P3.1.6}.
\sloppy

\ref{P3.1.7}.
Denote by $M$ the expression in the first representation.
Then, by \eqref{M}, we have $M\le{\Lambda}_U^\dag(\varphi_1,\varphi_2)$.
Conversely, let $V\in EI(U)$, $x_{1k}\in V$, $x_{2k}\in X$ for all $k\in\N$ and $d(x_{1k},x_{2k})\to0$ as $k\to+\infty$.
By \cref{L2.1}\,\ref{L2.1.6}, there exists a subset $V'\in EI(U)$ such that $V\in EI(V')$.
Hence, $x_{1k},x_{2k}\in V'$ for all large enough $k\in\N$.
Thus, $\Lambda_U^\dag(\varphi_1,\varphi_2)\leq \liminf_{k\to+\infty} (\varphi_1(x_{1k})+\varphi_2(x_{2k}))$, and
consequently, ${\Lambda}_U^\dag(\varphi_1,\varphi_2)\le M$.
This proves the first representation.
The other two representations follow from definition \eqref{M} and the first representation thanks to \cref{L2.1}\,\ref{L2.1.5}.
%
\end{proof}

\begin{remark}
\label{R3.1}
\begin{enumerate}
\item
Unlike \eqref{La} and \eqref{La0},
construction \eqref{M} is only meaningful when $\Int U\ne\es$; see \cref{P3.1}\,\ref{P3.1.5}.
\sloppy

\item
The restriction $x_1\in V$ in the first and second representations in part \ref{P3.1.7} of \cref{P3.1} can be replaced with $x_2\in V$.
Analogous replacements can be made in
the last two representations in part \ref{P3.1.9}.
\end{enumerate}
\end{remark}

The inequalities in parts \ref{P3.1.1} and \ref{P3.1.4} of \cref{P3.1} can be strict.
Inequality \eqref{La0qc}, opposite to the second inequality in part \ref{P3.1.1},
is an important qualification condition.
We are going to show that in some important situations it can be replaced by a weaker (thanks to the first inequality in part \ref{P3.1.4} and \cref{E3.3}) condition
\begin{gather}
\label{Mqc}
\inf_{U} (\varphi_1+\varphi_2)\le
{\Lambda}_U^\dag(\varphi_1,\varphi_2).
\end{gather}
Note that, unlike \eqref{La0qc}, inequality \eqref{Mqc} can be strict; see \cref{E3.2}.

The quantities compared in \eqref{La0qc} or \eqref{Mqc} are computed independently.
At the same time, it is important in some applications to ensure that, given an appropriate sequence of $(x_1,x_2)$ with $d(x_1,x_2)\to0$ as in \eqref{La0}, the corresponding $x$ approximating the infimum in the \LHS\ can be chosen close to $x_1$ and $x_2$ (which are close to each other because $d(x_1,x_2)\to0$).
To accommodate for this additional requirement, we are going to utilize the following definitions:
\begin{subequations}
\label{Th}
\begin{align}
\label{Th0}
{\Theta}_{U}^\circ(\varphi_1,\varphi_2):=&
\limsup_{\substack{d(x_1,x_2)\to0\\ x_1\in\dom\varphi_1\cap U,\,x_2\in\dom\varphi_2\cap U}}\; (\varphi_1\diamondsuit\varphi_2)_U(x_1,x_2),
\\
\label{Thdag}
{\Theta}_{U}^\dag(\varphi_1,\varphi_2):=&
\sup_{\substack{V\in EI(U)}}\; \limsup_{\substack{d(x_1,x_2)\to0\\ x_1\in\dom\varphi_1\cap V,\,x_2\in\dom\varphi_2\cap V}}\; (\varphi_1\diamondsuit\varphi_2)_U(x_1,x_2),
\end{align}
\end{subequations}
with the notation
\begin{align}
\label{phihat}
(\varphi_1\diamondsuit\varphi_2)_U(x_1,x_2):=
\inf_{x\in U} \max
\{d(x,x_1),d(x,x_2),
(\varphi_1+\varphi_2)(x)- \varphi_1(x_{1})-\varphi_2(x_{2})\}.
\end{align}
Thanks to assumption \eqref{nonempty} the latter expression is well defined as long as $x_1\in\dom\varphi_1$ and $x_2\in\dom\varphi_2$ (which is the case in \eqref{Th}).

In view of \eqref{phihat}, definitions \eqref{Th}
involve minimization of both $x\mapsto(\varphi_1+\varphi_2)(x)$ and $(x_1,x_2)\mapsto\varphi_1(x_1)+\varphi_2(x_2)$ as well as the distances $d(x,x_1)$ and $d(x,x_2)$.

\begin{remark}
\label{R3.2}
Since $d(x_1,x_2)\to0$ in \eqref{Th}, the definitions do not change if one of the terms $d(x,x_1)$ or $d(x,x_2)$ in definition \eqref{phihat} is dropped.
\end{remark}

Some elementary properties of the quantities \eqref{Th} and equivalent representations of the quantity \eqref{Thdag} are collected in the next proposition.

\begin{proposition}
\label{P3.2}
\begin{enumerate}
\item
\label{P3.2.1}
$0\le{\Theta}_{U}^\dag(\varphi_1,\varphi_2)\le {\Theta}_{U}^\circ(\varphi_1,\varphi_2)$.
\item
\label{P3.2.2}
If $\inf_{U} (\varphi_1+\varphi_2)>-\infty$ or ${\Lambda}_U^\circ(\varphi_1,\varphi_2)>-\infty$, then
\begin{align}
\label{P3.2-2a}
\inf_{U} (\varphi_1+\varphi_2)-
{\Lambda}_U^\circ(\varphi_1,\varphi_2) \le{\Theta}_{U}^\circ(\varphi_1,\varphi_2).
\end{align}
If $\inf_{U} (\varphi_1+\varphi_2)>-\infty$ or ${\Lambda}_U^\dag(\varphi_1,\varphi_2)>-\infty$, then
\begin{align}
\label{P3.2-2b}
\inf_{U} (\varphi_1+\varphi_2)-
{\Lambda}_U^\dag(\varphi_1,\varphi_2) \le{\Theta}_{U}^\dag(\varphi_1,\varphi_2).
\end{align}

\item
\label{P3.2.3}
The following estimates are true:
\begin{subequations}
\begin{align*}
{\Theta}_{U}^\circ(\varphi_1,\varphi_2)\ge&
\sup_{\substack{\al>0}}\; \limsup_{\substack{d(x_1,x_2)\to0;\, x_1,x_2\in U\\ \varphi_1(x_1)+\varphi_2(x_2)<\al}}\; (\varphi_1\diamondsuit\varphi_2)_U(x_1,x_2),
\\
{\Theta}_{U}^\dag(\varphi_1,\varphi_2)\ge&
\sup_{\substack{V\in EI(U)\\\al>0}}\; \limsup_{\substack{d(x_1,x_2)\to0;\, x_1,x_2\in V\\ \varphi_1(x_1)+\varphi_2(x_2)<\al}}\; (\varphi_1\diamondsuit\varphi_2)_U(x_1,x_2).
\end{align*}
\end{subequations}

\item
\label{P3.2.4}
The following representation holds true:
\begin{align*}
{\Theta}_{U}^\dag(\varphi_1,\varphi_2)&=
\sup\limits_{\substack{V\in EI(U)}}\; \limsup\limits_{\substack{d(x_1,x_2)\to0\\ x_1\in\dom\varphi_1\cap V,\,x_2\in\dom\varphi_2}}\; (\varphi_1\diamondsuit\varphi_2)_U(x_1,x_2).
\end{align*}
Moreover, $EI(U)$ in \eqref{Thdag} and the above representation can be replaced with $\EI(U)$.

\item
\label{P3.2.5}
If $X$ is a normed space, $\bx\in X$ and $\de>0$, then the following representations hold true:
\begin{align*}
{\Theta}_{B_\de(\bx)}^\dag(\varphi_1,\varphi_2)=
{\Theta}_{\overline B_\de(\bx)}^\dag(\varphi_1,\varphi_2)&=
\sup_{\substack{\rho\in(0,\de)}}\; \limsup_{\substack{d(x_1,x_2)\to0\\ x_1\in\dom\varphi_1\cap  B_\rho(\bx)\\x_2\in\dom\varphi_2\cap  B_\rho(\bx)}}\;
(\varphi_1\diamondsuit\varphi_2)_U(x_1,x_2)
\\&=
\sup_{\substack{\rho\in(0,\de)}}\; \limsup_{\substack{d(x_1,x_2)\to0\\ x_1\in\dom\varphi_1\cap  B_\rho(\bx),\,x_2\in\dom\varphi_2}}\;
(\varphi_1\diamondsuit\varphi_2)_U(x_1,x_2).
\end{align*}
Moreover, $B_\rho(\bx)$ in the above representations can be replaced with $\overline B_\rho(\bx)$.
\end{enumerate}
\end{proposition}

\begin{proof}
The assertions are direct consequences of definitions from \eqref{Th}.
For the first inequality in \ref{P3.2.1} in the case $\Int U=\es$, recall the convention $\sup\es_{\R_+}=0$.
For the inequalities in \ref{P3.2.2},
observe from \eqref{phihat} that
\begin{align*}
\inf_{U} (\varphi_1+\varphi_2)- \varphi_1(x_{1})-\varphi_2(x_{2})
\le
(\varphi_1\diamondsuit\varphi_2)_U(x_1,x_2)
\end{align*}
for all $x_1\in\dom\varphi_1$ and $x_2\in\dom\varphi_2$.
Assertion \ref{P3.2.3} clearly follows as $\varphi_1(x_1)+\varphi_2(x_2)<\alpha$ for some $\alpha>0$
immediately gives $x_1\in\dom\varphi_1$ and $x_2\in\dom\varphi_2$.
For the representations in \ref{P3.2.4} and \ref{P3.2.5},
reuse the arguments in the proof of \cref{P3.1}\,\ref{P3.1.7} and \ref{P3.1.9}.
\end{proof}

\begin{remark}
\label{R3.3}
\begin{enumerate}
\item
The assumptions in part \ref{P3.2.2} of \cref{P3.2} ensure that the \LHS s of \eqref{P3.2-2a} and \eqref{P3.2-2b} are well defined.
Recall that $\inf_{U} (\varphi_1+\varphi_2)<+\infty$ thanks to \eqref{nonempty}.
\item
\label{R3.3.2}
If $\inf_{U} (\varphi_1+\varphi_2)= {\Lambda}_U^\circ(\varphi_1,\varphi_2)=-\infty$ (or $\inf_{U} (\varphi_1+\varphi_2)= {\Lambda}_U^\dag(\varphi_1,\varphi_2)=-\infty$), then condition \eqref{La0qc} (or condition \eqref{Mqc}) is trivially satisfied.
\item
The restriction $x_1\in V$ in the representation in part \ref{P3.2.4} of \cref{P3.2} can be replaced with $x_2\in V$.
A similar replacement can be made in the second representation in part \ref{P3.2.5}.
\end{enumerate}
\end{remark}

Employing \eqref{Th}, the analogues of qualification conditions \eqref{La0qc} and \eqref{Mqc} can be written
as ${\Theta}_{U}^\circ(\varphi_1,\varphi_2)=0$ and ${\Theta}_{U}^\dag(\varphi_1,\varphi_2)=0$, respectively.
The next proposition examines these conditions.

\begin{proposition}
\label{P3.3}
\begin{enumerate}
\item
\label{P3.3.1}
${\Theta}_{U}^\dag(\varphi_1,\varphi_2)=0$
if and only if
\begin{equation*}
\sup_{\substack{V\in EI(U)}}\; \limsup_{\substack{d(x_1,x_2)\to0\\ x_1\in\dom\varphi_1\cap V,\,x_2\in\dom\varphi_2\cap V}}\; (\varphi_1\diamondsuit\varphi_2)_X(x_1,x_2)=0.
\end{equation*}

\item
\label{P3.3.2}
If ${\Theta}_{U}^\circ(\varphi_1,\varphi_2)=0$, then \eqref{La0qc} is satisfied.

\item
\label{P3.3.3}
If ${\Theta}_{U}^\dag(\varphi_1,\varphi_2)=0$, then \eqref{Mqc} is satisfied, and
${\Theta}_{U'}^\dag(\varphi_1,\varphi_2)=0$ for any $U'\subset U$.

\item
\label{P3.3.4}
Let $\bx\in X$.
The following conditions are equivalent:
\begin{enumerate}
\item
\label{P3.3.4a}
${\Theta}_{B_\de(\bx)}^\dag(\varphi_1,\varphi_2)=0$ for all sufficiently small $\de>0$;
\item
\label{P3.3.4c}
${\Theta}_{B_\de(\bx)}^\dag(\varphi_1,\varphi_2)=0$ for some $\de>0$.
\end{enumerate}
\end{enumerate}
\end{proposition}

\begin{proof}
\begin{enumerate}
\item
Set
\begin{equation*}
{\Theta}_{U,X}^\dag(\varphi_1,\varphi_2):=
\sup_{\substack{V\in EI(U)}}\; \limsup_{\substack{d(x_1,x_2)\to0\\ x_1\in\dom\varphi_1\cap V,\,x_2\in\dom\varphi_2\cap V}}\; (\varphi_1\diamondsuit\varphi_2)_X(x_1,x_2).
\end{equation*}
Since ${\Theta}_{U,X}^\dag(\varphi_1,\varphi_2)\le {\Theta}_{U}^\dag(\varphi_1,\varphi_2)$,
we only need to show that ${\Theta}_{U,X}^\dag(\varphi_1,\varphi_2)=0$
implies ${\Theta}_{U}^\dag(\varphi_1,\varphi_2)=0$.
Let ${\Theta}_{U,X}^\dag(\varphi_1,\varphi_2)=0$, $V\in EI(U)$, $x_{1k}\in\dom\varphi_1\cap V$,
$x_{2k}\in\dom\varphi_2\cap V$ for all $k\in\N$ and $d(x_{1k},x_{2k})\to 0$ as well as
$(\varphi_1\diamondsuit\varphi_2)_X(x_{1k},x_{2k})\to 0$ as $k\to+\infty$.
Then there exists a sequence $\{x_k\}\subset X$ such that conditions \eqref{ULC} are satisfied.
Since $x_{1k}\in V$, it follows from \cref{D2.1} that $x_{k}\in U$ for all sufficiently large $k\in\N$,
and consequently, $\lim_{k\to\infty} (\varphi_1\diamondsuit\varphi_2)_U(x_{1k},x_{2k})=0$.
Thus, ${\Theta}_{U}^\dag(\varphi_1,\varphi_2)=0$.

\item
This is a simple consequence of \cref{P3.2}\,\ref{P3.2.2} and \cref{R3.3}\,\ref{R3.3.2}.

\item
Let ${\Theta}_{U}^\dag(\varphi_1,\varphi_2)=0$.
Condition \eqref{Mqc} is a consequence of \cref{P3.2}\,\ref{P3.2.2} and \cref{R3.3}\,\ref{R3.3.2}.
Thanks to \ref{P3.3.1}, we have ${\Theta}_{U,X}^\dag(\varphi_1,\varphi_2)=0$.
Let $U'\subset U$.
Then ${\Theta}_{U',X}^\dag(\varphi_1,\varphi_2)\le {\Theta}_{U,X}^\dag(\varphi_1,\varphi_2)$.
Hence, ${\Theta}_{U',X}^\dag(\varphi_1,\varphi_2)=0$.
In view of \ref{P3.3.1}, the latter condition is equivalent to ${\Theta}_{U'}^\dag(\varphi_1,\varphi_2)=0$.

\item
Implication
\ref{P3.3.4a} \folgt \ref{P3.3.4c} is straightforward, while the converse implication follows from the second assertion in \ref{P3.3.3}.
\end{enumerate}
\end{proof}

\begin{remark}
\label{R3.4}
\begin{enumerate}
\item
\label{R3.4.1}
The analogues ${\Theta}_{U}^\circ(\varphi_1,\varphi_2)=0$ and ${\Theta}_{U}^\dag(\varphi_1,\varphi_2)=0$
of inequalities \eqref{La0qc} and \eqref{Mqc} are formulated as equalities.
This is because of the presence of the terms $d(x,x_1)$ and $d(x,x_2)$ in definition \eqref{phihat},
preventing $(\varphi_1\diamondsuit\varphi_2)_U(x_1,x_2)$ from being negative.
Besides, they do not allow to separate in \eqref{phihat} the terms containing $x_1,x_2$
on one hand and $x$ on the other hand (as in \eqref{La0qc} and~\eqref{Mqc}).
\item
Thanks to part \ref{P3.3.1} of \cref{P3.3} together with parts \ref{L2.1.5} and~\ref{L2.1.8} of \cref{L2.1},
$B_\de(\bx)$ in part \ref{P3.3.4} of \cref{P3.3} can be replaced with $\overline B_\de(\bx)$.
\end{enumerate}
\end{remark}

The examples below illustrate the computation of the `decoupling quantities' \eqref{La}, \eqref{La0}, \eqref{M} and~\eqref{Th}.

\begin{example}
\label{E3.1}
Let \lsc\ convex functions $\varphi_1,\varphi_2\colon\R^2\to\R_\infty$ be given by
\begin{align}
\label{E3.1-1}
\forall (x,y)\in\R^2\colon\quad
\varphi_1(x,y):=
\begin{cases}
-x&\text{if }y\ge x^2,
\\
+\infty&\text{otherwise},
\end{cases}
\quad
\varphi_2(x,y):=
\begin{cases}
0&\text{if }y\le0,
\\
+\infty&\text{otherwise}.
\end{cases}
\end{align}
Then $\dom\varphi_1\cap\dom\varphi_2=\{(0,0)\}$, and $\varphi_1(0,0)=\varphi_2(0,0)=0$.
Let $U\subset\R^2$ and $(0,0)\in\Int U$.
Then $\inf_{U}(\varphi_1+\varphi_2)= (\varphi_1+\varphi_2)(0,0)=0$.
If $(x_1,y_1)\in\dom\varphi_1$, $(x_2,y_2)\in\dom\varphi_2$,
and
$\norm{(x_1,y_1)-(x_2,y_2)}\to0$, then $(x_1,y_1)\to(0,0)$ and $(x_2,y_2)\to(0,0)$.
Hence,
$\Lambda_U(\varphi_1,\varphi_2)= \Lambda_U^\circ(\varphi_1,\varphi_2)=
{\Lambda}_{U}^\dag(\varphi_1,\varphi_2)= {\Theta}_{U}^\circ(\varphi_1,\varphi_2)=
{\Theta}_{U}^\dag(\varphi_1,\varphi_2)=0$.
\end{example}

The next example shows that inequality \eqref{Mqc} and, thus,
the first inequality in \cref{P3.1}\,\ref{P3.1.4}
can be strict.

\begin{example}
\label{E3.2}
Let \lsc\ convex functions $\varphi_1,\varphi_2\colon\R\to\R$ be given by
\begin{align*}
\forall x\in\R\colon\quad
\varphi_1(x):=
\begin{cases}
0&\text{if }x\le0,
\\
1&\text{otherwise},
\end{cases}
\quad
\varphi_2(x):=0,
\end{align*}
and $U:=[0,1]$.
Thus, $\varphi_1+\varphi_2=\varphi_1$ and
$\inf_{U}(\varphi_1+\varphi_2)=
\Lambda_U(\varphi_1,\varphi_2)= \Lambda_U^\circ(\varphi_1,\varphi_2)=0$, while
${\Lambda}_{U}^\dag(\varphi_1,\varphi_2)=1$.
If $x_1,x_2\in U$, then
\begin{align*}
(\varphi_1\diamondsuit\varphi_2)_U(x_1,x_2)
&=
\inf_{x\in U} \max\{|x-x_1|,|x-x_2|,\varphi_1(x)- \varphi_1(x_{1})\}\\
&=
\begin{cases}
	x_2-x_1	&\text{if }x_1=0,
	\\
	|x_2-x_1|/2	&\text{otherwise},
\end{cases}
\end{align*}
and consequently, $(\varphi_1\diamondsuit\varphi_2)_U(x_1,x_2)\to0$ as $|x_1-x_2|\to0$.
It follows that ${\Theta}_{U}^\circ(\varphi_1,\varphi_2)=
{\Theta}_{U}^\dag(\varphi_1,\varphi_2)=0$.
Particularly, \eqref{P3.2-2a} holds as equality while inequality \eqref{P3.2-2b} is strict.
\sloppy

\end{example}

The next example illustrates a situation where, again, the
first inequality in \cref{P3.1}\,\ref{P3.1.4} is strict
while both \eqref{P3.2-2a} and \eqref{P3.2-2b} hold as equalities.

\begin{example}
\label{E3.3}
Let $\de>0$, and \lsc\ functions $\varphi_1,\varphi_2\colon\R\to\R_\infty$ be given by
\begin{align*}
\forall x\in\R\colon\quad
\varphi_1(x):=
\begin{cases}
\frac{\de}{\de-x}&\text{if }x<\de,
\\
+\infty&\text{otherwise},
\end{cases}
\quad
\varphi_2(x):=
\begin{cases}
\frac{\de}{x-\de}&\text{if }x<\de,
\\
+\infty&\text{otherwise}.
\end{cases}
\end{align*}
Then
\begin{align*}
\forall x\in\R\colon\quad
(\varphi_1+\varphi_2)(x):=
\begin{cases}
0&\text{if }x<\de,
\\
+\infty&\text{otherwise}.
\end{cases}
\end{align*}
Consider numbers $\eta$ and $\rho$ satisfying
$0<\eta<\rho\le\de$.
A direct calculation gives
\begin{align*}
\la_{\rho,\eta}:=\inf_{\substack{x_1,x_2\in(-\rho,\rho),\, |x_1-x_2|<\eta}} (\varphi_1(x_1)+\varphi_2(x_2))=&
\lim_{x_1\downarrow(\rho-\eta)}\varphi_1(x_1)+ \lim_{x_2\uparrow\rho}\varphi_2(x_2)
\\=&
\begin{cases}
-\frac{\de\eta}{(\de-\rho+\eta)(\de-\rho)}
&\text{if }\rho<\de,
\\
-\infty&\text{if }\rho=\de.
\end{cases}
\end{align*}
Let $U:=(-\delta,\delta)$. Then
$\inf_{U}(\varphi_1+\varphi_2)=0$, and
\[
\begin{aligned}
\Lambda_U(\varphi_1,\varphi_2)&= \Lambda_U^\circ(\varphi_1,\varphi_2)= \lim_{\eta\downarrow0} \la_{\de,\eta}=-\infty,&\quad
{\Lambda}_{U}^\dag(\varphi_1,\varphi_2)&=
\inf_{\rho\in(0,\de)}\; \lim_{\eta\downarrow0} \la_{\rho,\eta}=0,&
\\
{\Theta}_{U}^\circ(\varphi_1,\varphi_2)&= -\lim_{\eta\downarrow0} \la_{\de,\eta}=+\infty,& \quad
{\Theta}_{U}^\dag(\varphi_1,\varphi_2)&= -\inf_{\rho\in(0,\de)}\; \lim_{\eta\downarrow0} \la_{\rho,\eta}=0.&
\end{aligned}
\]
Thus, $\Lambda_U^\circ(\varphi_1,\varphi_2)< {\Lambda}_{U}^\dag(\varphi_1,\varphi_2)$.
At the same time, both inequalities in \cref{P3.2}\,\ref{P3.2.2} hold as equalities.
\end{example}

The next example shows that the inequalities in parts \ref{P3.2.2} and \ref{P3.2.3} of \cref{P3.2} can be strict.

\begin{example}
\label{E3.4}
Let \lsc\
functions $\varphi_1,\varphi_2\colon\R^2\to\R_\infty$ be given by
	\[
		\forall (x,y)\in\R^2\colon\quad
		\varphi_1(x,y)
		:=
		\begin{cases}
			0	&\text{if }	x=y=0,\\
			1/x	&\text{if }	x>0,\\
			+\infty	&	\text{otherwise,}
		\end{cases}
		\qquad
		\varphi_2(x,y):=\varphi_1(-x,y).
	\]
Then $\varphi_1(x,y)\ge0$ and $\varphi_2(x,y)\ge0$ for all $(x,y)\in\R^2$,
$\dom\varphi_1\cap\dom\varphi_2=\{(0,0)\}$, $\varphi_1(0,0)=\varphi_2(0,0)=0$.
Let $U\subset\R^2$ and $(0,0)\in\Int U$.
Then $\inf_{U}(\varphi_1+\varphi_2)= (\varphi_1+\varphi_2)(0,0)=0$, and
$\Lambda_U(\varphi_1,\varphi_2)= \Lambda_U^\circ(\varphi_1,\varphi_2)=
{\Lambda}_{U}^\dag(\varphi_1,\varphi_2)=0$.

Let $U:=\delta\B$ for some $\delta>0$.
Choose a $\rho\in(0,\de)$, and set $V:=\rho\overline\B$ and $x_k:=\rho/k$ for all $k\in\N$.
Then $V\in EI(U)$, $(x_{k},\rho),(-x_{k},\rho)\in V$, $\varphi_1(x_{k},\rho)=\varphi_2(-x_{k},\rho)=1/x_k>0$ for all $k\in\N$, and
$\norm{(x_{k},\rho)-(-x_{k},\rho)}=2x_k\to 0$ as $k\to+\infty$.
Thus, $(\varphi_1\diamondsuit\varphi_2)_U ((x_{k},\rho),(-x_{k},\rho))=\max\{ \norm{(x_{k},\rho)},-2/x_k\}= \norm{(x_{k},\rho)}\ge\rho$
for all $k\in\N$.
Hence, ${\Theta}_{U}^\circ(\varphi_1,\varphi_2)\ge {\Theta}_{U}^\dag(\varphi_1,\varphi_2)\ge\de
$.
Note that the \RHS s of the inequalities in \cref{P3.2}\,\ref{P3.2.3} equal $0$
since $\varphi_1(x_{k},\rho)+\varphi_2(-x_{k},\rho)= 2/x_k=2k/\rho$, and this quantity becomes greater than any given number $\al>0$ when $k$ is large enough.

\end{example}

The final example of this section shows that
if condition \eqref{nonempty} is violated, then 
both sides of the inequalities \eqref{La0qc} and \eqref{Mqc}
may equal $+\infty$.
It also provides a situation when the inequalities in \cref{P3.2}\,\ref{P3.2.3} are strict.

\begin{example}
Let \lsc\ convex functions $\varphi_1,\varphi_2\colon\R\to\R_\infty$ be given by
\[
\forall x\in\R\colon\quad
\varphi_1(x):=
\begin{cases}
0&\text{if }x\le0,\\
+\infty&\text{otherwise,}
\end{cases}
\qquad
\varphi_2(x):=
\begin{cases}
1/x&\text{if }x>0,\\
+\infty&\text{otherwise,}
\end{cases}
\]
and pick a subset $U\subset\R$ with $0\in\Int U$.
Then condition \eqref{nonempty} fails, and consequently, $\inf_{U}(\varphi_1+\varphi_2)=+\infty$.
If $x_1\in\dom\varphi_1$, $x_2\in\dom\varphi_2$,
and
$|x_1-x_2|\to0$, then $x_2\downarrow0$, and consequently, $\varphi_1(x_1)=0$ and $\varphi_2(x_2)\to+\infty$.
Hence,
${\Lambda}_U^\circ(\varphi_1,\varphi_2)= {\Lambda}_U^\dag(\varphi_1,\varphi_2)=+\infty$.
Moreover,
$(\varphi_1\diamondsuit\varphi_2)_U(x_1,x_2)=+\infty$ for all $x_1,x_2\in\R$, and consequently,
${\Theta}_{U}^\circ(\varphi_1,\varphi_2)= {\Theta}_{U}^\dag(\varphi_1,\varphi_2)=+\infty$.
At the same time,
for any $\al>0$, we have $\varphi_2(x_2)>\al$ when $x_2>0$ is small enough; hence, by convention,
the \RHS s of the inequalities in \cref{P3.2}\,\ref{P3.2.3} equal $0$.
\end{example}

\section{Uniform and quasiuniform lower semicontinuity}
\label{S4}

In this section we discuss certain uniform and quasiuniform lower semicontinuity properties
of pairs of functions resulting from definitions \eqref{La0}, \eqref{M} and \eqref{Th} of uniform infima which turn out to be
beneficial in the context of the fuzzy multiplier and calculus rules.

As in \cref{S3}, we consider a pair of functions
$\varphi_1,\varphi_2\colon X\to\R_\infty$ on a metric space $X$ and a subset $U\subset X$ satisfying condition \eqref{nonempty}.

\begin{definition}
\label{D4.1}
The pair $(\varphi_1,\varphi_2)$ is \begin{enumerate}
\item
\label{D4.1.1}
\emph{uniformly lower semicontinuous}
on $U$ if condition \eqref{La0qc} is satisfied;

\item
\label{D4.1.2}
\emph{quasiuniformly lower semicontinuous}
on $U$ if condition \eqref{Mqc} is satisfied;

\item
\label{D4.1.3}
\emph{firmly uniformly lower semicontinuous}
on $U$ if ${\Theta}_{U}^\circ(\varphi_1,\varphi_2)=0$;

\item
\label{D4.1.4}
\emph{firmly quasiuniformly lower semicontinuous}
on $U$ if ${\Theta}_{U}^\dag(\varphi_1,\varphi_2)=0$;

\item
\label{D4.1.5}
uniformly/quasiuniformly/firmly uniformly/firmly quasiuniformly lower semicontinuous
near a point $\bx\in\dom\varphi_1\cap\dom\varphi_2$ if it is uniformly/quasiuniformly/firmly uniform\-ly/firmly quasiuniformly lower semicontinuous on $\overline B_\de(\bx)$
for all sufficiently small $\de>0$.
\end{enumerate}
\end{definition}

\begin{remark}
\label{R4.1}
\begin{enumerate}
\item
\cref{D4.1}\,\ref{D4.1.1} follows \cite[Definition~2.6]{BorweinZhu96} (see also
\cite[Remark~2]{BorweinIoffe1996} and \cite[Section~2.3]{Las01}).
Some sufficient conditions for this property can be found in 
\cite{Las01} and \cite[Section~1.6.4]{Pen13}.

\item
\label{R4.1.2}
Thanks to \cref{P3.1}\,\ref{P3.1.4}, \cref{P3.2}\,\ref{P3.2.1} and \cref{P3.3}\,\ref{P3.3.2} and \ref{P3.3.3},
the properties in \cref{D4.1}\,\ref{D4.1.1}--\ref{D4.1.4} are related as follows:
\begin{align}
\label{R4.1-2}
\ref{D4.1.3} \Rightarrow \ref{D4.1.4} \Rightarrow \ref{D4.1.2} \AND
\ref{D4.1.3} \Rightarrow \ref{D4.1.1} \Rightarrow \ref{D4.1.2}.
\end{align}
The implications can be strict; see \cref{E3.3,E3.4}.

\item
In view of definitions \eqref{La0}, \eqref{M}, \eqref{Th} and \eqref{phihat}
the key conditions \eqref{La0qc}, \eqref{Mqc}, ${\Theta}_{U}^\circ(\varphi_1,\varphi_2)=0$ and ${\Theta}_{U}^\dag(\varphi_1,\varphi_2)=0$ in \cref{D4.1}
compare (some form of) the infimum of the sum $x\mapsto(\varphi_1+\varphi_2)(x)$
with (some form of) the lower limit of the decoupled sum $(x_1,x_2)\mapsto\varphi_1(x_1)+\varphi_2(x_2)$ as $d(x_1,x_2)\to0$.
These conditions prevent
the lower limit of the decoupled sum from being smaller than
the infimum of the sum.

\item
The properties in parts \ref{D4.1.1} and \ref{D4.1.2} of \cref{D4.1}
are defined via inequalities \eqref{La0qc} and \eqref{Mqc},
while the definition of their firm counterparts in parts \ref{D4.1.3} and \ref{D4.1.4}
use equalities ${\Theta}_{U}^\circ(\varphi_1,\varphi_2)=0$ and
${\Theta}_{U}^\dag(\varphi_1,\varphi_2)=0$.
This is because of the presence of the terms $d(x,x_1)$ and $d(x,x_2)$ in definition \eqref{phihat}
which prevent $(\varphi_1\diamondsuit\varphi_2)_U(x_1,x_2)$ from being negative; see \cref{R3.4}\,\ref{R3.4.1}.

\item
\label{R4.1.5}
Slightly weaker versions of the `firm' properties can be considered, corresponding to replacing
${\Theta}_{U}^\circ(\varphi_1,\varphi_2)$ and
${\Theta}_{U}^\dag(\varphi_1,\varphi_2)$ in parts \ref{D4.1.3} and \ref{D4.1.4} of \cref{D4.1}
by the \RHS s of the respective inequalities in \cref{P3.2}\,\ref{P3.2.3}.
\item
\cref{D4.1} together with definitions \eqref{La}, \eqref{La0}, \eqref{M}, \eqref{Th} and \eqref{phihat}
can be easily extended from pairs to arbitrary finite families of functions.
\end{enumerate}
\end{remark}

The characterizations of the uniform and quasiuniform lower semicontinuity properties
in the next proposition are consequences of \cref{D4.1}, definitions \eqref{La0}, \eqref{M}, \eqref{Th}, \eqref{phihat},
\cref{P3.1}\,\ref{P3.1.7}, \cref{P3.2}\,\ref{P3.2.4}, \cref{P3.3}\,\ref{P3.3.1} and \cref{R3.2}.

\begin{proposition}
\label{P4.1}
The pair $(\varphi_1,\varphi_2)$ is \begin{enumerate}
\item
uniformly lower semicontinuous
on $U$ if and only if, for any
$\eps>0$, there exists an $\eta>0$ such that, for any
$x_1,x_2\in U$
with $d(x_1,x_2)<\eta$, there is an $x\in U$ satisfying
\begin{gather}
\label{P4.1-1}
(\varphi_1+\varphi_2)(x)< \varphi_1(x_{1})+\varphi_2(x_{2})+\eps;
\end{gather}

\item
\label{P4.1.2}
quasiuniformly lower semicontinuous
on $U$ if and only if, for any
$V\in EI(U)$
and $\eps>0$, there exists an $\eta>0$ such that, for any
$x_1\in V$ and $x_2\in X$
with $d(x_1,x_2)<\eta$, there is an $x\in U$ satisfying condition \eqref{P4.1-1};

\item
\label{P4.1.3}
{firmly uniformly lower semicontinuous}
on $U$ if and only if, for any $\eps>0$, there exists an $\eta>0$ such that, for any
$x_1\in\dom\varphi_1\cap U$ and $x_2\in\dom\varphi_2\cap U$
with $d(x_1,x_2)<\eta$, there is an $x\in U\cap B_\eps(x_1)$ satisfying condition \eqref{P4.1-1};

\item
\label{P4.1.4}
{firmly quasiuniformly lower semicontinuous}
on $U$ if and only if, for any $V\in EI(U)$ and $\eps>0$, there exists an $\eta>0$ such that, for any
$x_1\in\dom\varphi_1\cap V$ and $x_2\in\dom\varphi_2$
with $d(x_1,x_2)<\eta$, there is an $x\in B_\eps(x_1)$ satisfying condition \eqref{P4.1-1}.
\end{enumerate}
\end{proposition}

Thanks to \cref{P3.3}\,\ref{P3.3.4}, it is possible to state a simplified characterization of
firm quasiuniform lower semicontinuity near a point.
It replaces the closed ball $\overline B_\de(\bar x)$ with the open ball $B_\de(\bar x)$ and requires condition $\Theta^\dag_{B_\de(\bar x)}(\varphi_1,\varphi_2)=0$ to hold not for all, but for \emph{some} $\de>0$.

\begin{proposition}
\label{prop:firm_quasiuniform_lsc_near_point}
	The pair $(\varphi_1,\varphi_2)$ is
	firmly quasiuniformly lower semicontinuous
	near $\bx\in\dom\varphi_1\cap\dom\varphi_2$ if and only if $\Theta^\dag_{B_\de(\bar x)}(\varphi_1,\varphi_2)=0$ for some $\de>0$.
\end{proposition}


The next proposition gives sequential reformulations of the characterizations of quasiuniform lower semicontinuity properties from \cref{P4.1}.

\begin{proposition}
\label{P4.3}
The pair $(\varphi_1,\varphi_2)$ is \begin{enumerate}
\item
\label{P4.3.1}
uniformly lower semicontinuous
on $U$ if and only if, for any sequences
$\{x_{1k}\}\subset\dom\varphi_1\cap U$ and $\{x_{2k}\}\subset\dom\varphi_2\cap U$
satisfying $d(x_{1k},x_{2k})\to0$ as $k\to+\infty$, there exists a sequence $\{x_k\}\subset U$ such that condition \eqref{ULCb} is satisfied;

\item
\label{P4.3.2}
quasiuniformly lower semicontinuous
on $U$ if and only if, for any sequences
$\{x_{1k}\}\subset\dom\varphi_1\cap U$ and $\{x_{2k}\}\subset\dom\varphi_2$ satisfying $\{x_{1k}\}\in EI(U)$ and $d(x_{1k},x_{2k})\to0$ as $k\to+\infty$, there exists a sequence $\{x_k\}\subset U$ such that
condition \eqref{ULCb} is satisfied;

\item
\label{P4.3.3}
firmly uniformly lower semicontinuous
on $U$ if and only if, for any sequences
$\{x_{1k}\}\subset\dom\varphi_1\cap U$ and
$\{x_{2k}\}\subset\dom\varphi_2\cap U$
satisfying $d(x_{1k},x_{2k})\to0$ as $k\to+\infty$, there exists a sequence $\{x_k\}\subset U$ such that
conditions \eqref{ULC} are satisfied;

\item
\label{P4.3.4}
firmly quasiuniformly lower semicontinuous
on $U$ if and only if, for any sequences
$\{x_{1k}\}\subset\dom\varphi_1\cap U$ and
$\{x_{2k}\}\subset\dom\varphi_2$
satisfying $\{x_{1k}\}\in EI(U)$ and $d(x_{1k},x_{2k})\to0$ as $k\to+\infty$, there exists a sequence $\{x_k\}\subset X$ such that
conditions \eqref{ULC} are satisfied.
\end{enumerate}
\end{proposition}

\begin{remark}
\label{R4.2}
\begin{enumerate}
\item
\label{R4.2.1}
Unlike the characterizations in the other parts of \cref{P4.1,P4.3},
those in parts \ref{P4.1.4} of these statements do not require explicitly that $x$ or $x_k$ belong to $U$ or $B_\de(\bx)$.
However, in view of \cref{P3.3}\,\ref{P3.3.1}, the mentioned conditions are automatically satisfied in these characterizations.

\item
Condition \eqref{P4.1-1} automatically implies that $x\in\dom\varphi_1\cap\dom\varphi_2$.
The condition obviously only needs to be checked for $x_1\in\dom\varphi_1$ and $x_2\in\dom\varphi_2$.
In view of \cref{P4.1}, the properties in parts \ref{D4.1.2} and \ref{D4.1.4} of \cref{D4.1} can only be meaningful when $\dom\varphi_1\cap V\ne\es$ and $\dom\varphi_2\cap V\ne\es$ for some $V\in EI(U)$.

\item
\label{R4.2.3}
Due to condition $x\in B_\eps(x_1)$ involved in \cref{P4.1}\,\ref{P4.1.4} and \cref{prop:firm_quasiuniform_lsc_near_point},
it looks as if the point $x_1$ plays a special role in the firm quasiuniform lower semicontinuity property.
In fact, both $x_1$ and $x_2$ contribute equally to this property (see definition \eqref{phihat} and \cref{R3.2}), and the mentioned condition can be replaced there with
$x\in B_\eps(x_1)\cap B_\eps(x_2)$.

\item
Another example of a seeming lack of symmetry appears in parts \ref{P4.1.2} and \ref{P4.1.4} of \cref{P4.1}:
they require $x_1$ to belong to some essentially interior subset $V$, while there are no such restrictions on $x_2$.
This is because they use equivalent representations of ${\Theta}_{U}^\dag(\varphi_1,\varphi_2)$ from \cref{P3.2}\,\ref{P3.2.4}.
Recall that definition \eqref{Thdag} of ${\Theta}_{U}^\dag(\varphi_1,\varphi_2)$ requires both $x_1$ and $x_2$ to belong to $V$.
Thus, these characterizations can be rewritten in a (formally slightly more restrictive) symmetric form.
A similar observation applies to the characterization in \cref{prop:firm_quasiuniform_lsc_near_point}.

\item
\label{R4.2.5}
If $X$ is a normed space, and $U=B_\de(\bx)$ or $U=\overline B_\de(\bx)$ for some $\bx\in X$ and $\de>0$, one can use the characterizations from \cref{P3.1}\,\ref{P3.1.9}
to replace the arbitrary essentially interior subsets $V$
in part \ref{P4.1.2} of \cref{P4.1}
by either the family of closed balls $\overline B_\rho(\bar x)$ or the family of open balls $B_\rho(\bar x)$
with $\rho\in(0,\de)$.

\item
In view of \cref{L2.1}\,\ref{L2.1.5} collection $EI(U)$ in parts \ref{P4.1.2} and \ref{P4.1.4} of \cref{P4.1} can be replaced with its sub-collection $\EI(U)$.

\item
\label{R4.2.6}
In view of the characterization in \cref{P4.3}\,\ref{P4.3.3}
the \emph{(ULC) property} \cite[Definition~6]{BorweinIoffe1996} 
(also known as \emph{sequential uniform lower semicontinuity}
\cite[Definition~3.3.17]{BorZhu05}) at $\bx$ is equivalent to the firm uniform lower semicontinuity on $\overline B_\de(\bx)$ for some $\de>0$.
Thanks to parts \ref{P4.3.1} and \ref{P4.3.3} of \cref{P4.3}, a pair of functions is uniformly lower semicontinuous
(resp., firmly uniformly lower semicontinuous) on $X$ if it is quasicoherent (resp., coherent) \cite[Lemma~1.124]{Pen13}.

\item
\label{R4.2.7}
The properties in \cref{D4.1} are rather weak.
This is illustrated by the examples in \cref{S3}.
For pairs of functions $\varphi_1$ and $\varphi_2$ and corresponding sets $U$ 
in \cref{E3.1,E3.2,E3.3,E3.4}, it has been shown that $\inf_{U}(\varphi_1+\varphi_2)\le {\Lambda}_{U}^\dag(\varphi_1,\varphi_2)$,
i.e., all the pairs are quasiuniformly lower semicontinuous on the respective sets $U$.
Moreover, in \cref{E3.1,E3.2,E3.3}, it also holds
${\Theta}_{U}^\dag(\varphi_1,\varphi_2)=0$, i.e., the pairs are actually firmly quasiuniformly lower semicontinuous on $U$,
while in \cref{E3.1,E3.2}, ${\Theta}_{U}^\circ(\varphi_1,\varphi_2)=0$, i.e., the pairs are uniformly lower semicontinuous on $U$.
At the same time,
in \cref{E3.4},
${\Theta}_{U}^\dag(\varphi_1,\varphi_2)>0$, thus, illustrating that firm quasiuniform lower semicontinuity is indeed a stronger property than its non-firm counterpart.

Furthermore, in \cref{E3.3}, we have
$\Lambda_U^\circ(\varphi_1,\varphi_2)= -\infty$
while $\inf_{U}(\varphi_1+\varphi_2)=0$, i.e., condition \eqref{La0qc} is violated, and the functions are not uniformly lower semicontinuous
(see \cref{R4.1}\,\ref{R4.1.2}).
As a consequence,
the assertion of
\cite[Exercise~3.3.8]{BorZhu05} is not correct (see item \ref{R4.2.6} above).
\end{enumerate}
\end{remark}

We now show that firm uniform and firm quasiuniform lower semicontinuity properties are stable under uniformly continuous perturbations of the involved functions.

\begin{proposition}
\label{P4.4}
Suppose that $(\varphi_1,\varphi_2)$ is firmly uniformly (resp., firmly quasiuniformly) lower semicontinuous and $g\colon X\to\R$ is uniformly continuous on $U$.
Then $(\varphi_1,\varphi_2+g)$ is firmly uniformly (resp., firmly quasiuniformly) lower semicontinuous on $U$.
\end{proposition}

\begin{proof}
We employ
\cref{P4.3}\,\ref{P4.3.3} and \ref{P4.3.4}.
Let sequences
$\{x_{1k}\}\subset\dom\varphi_1$ and
$\{x_{2k}\}\subset\dom(\varphi_2+g)$
satisfy
$\{x_{1k}\},\{x_{2k}\}\subset U$
(resp., $\{x_{1k}\}\in EI(U)$) and $d(x_{1k},x_{2k})\to0$ as $k\to+\infty$.
Then $\{x_{2k}\}\subset\dom\varphi_2$.
Hence, there exists a sequence $\{x_k\}\subset U$ such that
conditions \eqref{ULC} are satisfied.
The uniform continuity of $g$ on $U$ gives
$g(x_k)-g(x_{2k})\to0$,
and in view of \eqref{ULCb}, we have
\begin{align*}
\limsup_{k\to +\infty} \big((\varphi_1+\varphi_2+g)(x_k)&- \varphi_1(x_{1k})-(\varphi_2+g)(x_{2k})\big)
\\
&=
\limsup_{k\to +\infty} \big((\varphi_1+\varphi_2)(x_k)- \varphi_1(x_{1k})-\varphi_2(x_{2k})\big)\le0.
\end{align*}
i.e., $(\varphi_1,\varphi_2+g)$ is firmly uniformly (resp., firmly quasiuniformly) lower semicontinuous on~$U$.
\end{proof}

Let us point out that the proof of \cref{P4.4} heavily exploits the nature of firm uniform/quasiuniform lower semicontinuity.
The assertion may not be true if it is replaced by its non-firm counterpart.

The next proposition collects several sufficient conditions for firm uniform lower semicontinuity.

\begin{proposition}
\label{P4.5}
The pair $(\varphi_1,\varphi_2)$ is {firmly uniformly lower semicontinuous} on $U$ provided that one of the following conditions is satisfied:
\begin{enumerate}
\item
\label{P4.5.1}
there is a $c\in\R$ such that
$\varphi_2(x)=c$ for all $x\in\dom\varphi_1\cap U$ and $\varphi_2(x)\ge c$ for all $x\in U\setminus\dom\varphi_1$
(or simply $\varphi_2(x)=c$ for all $x\in U$);
\item
\label{P4.5.2}
$\dom\varphi_2\cap U=\{\bar x\}$,
$\bar x\in\dom\varphi_1$
and
$\varphi_1$ is lower semicontinuous at $\bx$;
\item
\label{P4.5.3}
$\varphi_2$ is uniformly continuous on $U$.
\end{enumerate}
\end{proposition}

\begin{proof}
\begin{enumerate}
\item
The assertion is a consequence of \cref{P4.1}\,\ref{P4.1.3}.
Let $\varepsilon>0$.
Set $\eta:=\eps$.
Given any points $x_1\in\dom\varphi_1\cap U$ and $x_2\in\dom\varphi_2\cap U$
with $d(x_1,x_2)<\eta$, we have $\varphi_2(x_1)=c$,
$\varphi_2(x_2)\ge c$ and
condition \eqref{P4.1-1} is satisfied with $x:=x_{1}$.

\item
Let sequences $\{x_{1k}\}\subset\dom\varphi_1\cap U$ and $\{x_{2k}\}\subset\dom\varphi_2\cap U$
satisfy $d(x_{1k},x_{2k})\to0$ as $k\to+\infty$.
Then $x_{2k}=\bx$.
Thus, $x_{1k}\to\bar x$ as $k\to+\infty$ and, in view of the lower semicontinuity of $\varphi_1$,
the conditions in \cref{P4.3}\,\ref{P4.3.3} are satisfied with $x_k:=\bx$ for all $k\in\N$.
\item
This is a consequence of \ref{P4.5.1} and \cref{P4.4}.
\end{enumerate}
\end{proof}

\begin{remark}
\begin{enumerate}
\item
Condition \ref{P4.5.1} in \cref{P4.5} is satisfied if $\varphi_2$ is constant everywhere or is the indicator function of a set containing $\dom\varphi_1\cap U$.
Function $\varphi_1$ in condition \ref{P4.5.1}
does not have to be lower semicontinuous in the conventional sense.
\item
\cref{P4.5} with condition \ref{P4.5.3} strengthens \cite[Proposition~2.7.1]{BorweinZhu96}
and is similar to \cite[Proposition~2.1\,(d)]{Las01} and \cite[Exercise~3.2.9]{BorZhu05}
(which use \eqref{La} instead of \eqref{La0} in the qualification condition \eqref{La0qc}).
\end{enumerate}
\end{remark}

The next proposition
exploits compactness assumptions 
(often referred to as \emph{inf-compact\-ness})
in order to guarantee
the uniform lower semicontinuity properties in \cref{D4.1}.

\begin{proposition}
\label{P4.6}
Suppose that
$\varphi_1$ and $\varphi_2$ are \lsc\ on $U$ and $\inf_U\varphi_2>-\infty$.
\begin{enumerate}
\item
\label{P4.6.1}
The pair $(\varphi_1,\varphi_2)$ is {uniformly lower semicontinuous}
on $U$ if $\{x\in U\,|\,\varphi_1(x)\leq c\}$ is compact for each $c\in\R$.
If, additionally,
\begin{align}
\label{P4.6-3}
\limsup_{\substack{d(x_1,x_2)\to0\\ x_1\in\dom\varphi_1\cap U,\,x_2\in\dom\varphi_2\cap U}}\; (\varphi_1(x_1)+\varphi_2(x_2))<+\infty,
\end{align}
then it is {firmly uniformly lower semicontinuous} on $U$.
\item
\label{P4.6.2}
The pair $(\varphi_1,\varphi_2)$ is {quasiuniformly lower semicontinuous}
on $U$ if $\{x\in\cl V\,|\,\varphi_1(x)\leq c\}$ is compact for each $V\in EI(U)$ and $c\in\R$.
If, additionally,
\begin{align}
\label{P4.6-4}
\sup_{\substack{V\in EI(U)}}\; \limsup_{\substack{d(x_1,x_2)\to0\\ x_1\in\dom\varphi_1\cap V,\,x_2\in\dom\varphi_2}}\; (\varphi_1(x_1)+\varphi_2(x_2))<+\infty,
\end{align}
then it is {firmly quasiuniformly lower semicontinuous} on $U$.
\end{enumerate}
\end{proposition}

\begin{proof}
We consider sequences
$\{x_{1k}\}\subset\dom\varphi_1$ and
$\{x_{2k}\}\subset\dom\varphi_2$
such that
$d(x_{1k},x_{2k})\to0$ as $k\to+\infty$.
We assume that $\{x_{1k}\},\{x_{2k}\}\subset U$ in case \ref{P4.6.1},
and $\{x_{1k}\}\subset V$ for some $V\in EI(U)$
in case \ref{P4.6.2}.
The latter assumption implies that
$x_{2k}\in U$
for all sufficiently large $k\in\N$.
Thus, in all cases we have $\varphi_2(x_{2k})\ge m:=\inf_U\varphi_2$
for all sufficiently large $k\in\N$.
If $\varphi_1(x_{1k})+\varphi_2(x_{2k})\to+\infty$ as $k\to+\infty$, then, taking any $\hat x\in\dom\varphi_1\cap\dom\varphi_2\cap U$,
condition \eqref{ULCb} is satisfied with $x_k:=\hat x$ for all $k\in\N$.
Assume without loss of generality that $\al:=\limsup_{k\to\infty} (\varphi_1(x_{1k})+\varphi_2(x_{2k}))<+\infty$
(this is automatically ensured by the additional conditions \eqref{P4.6-3} and \eqref{P4.6-4}).
Then, for all sufficiently large $k\in\N$, we have
$\varphi_1(x_{1k})<c:=\al-m+1$, i.e., a tail of the sequence $\{x_{1k}\}$
belongs to $\{x\in U\,|\,\varphi_1(x)\leq c\}$.
Moreover, in case \ref{P4.6.2} it belongs to $\{x\in \cl V\,|\,\varphi_1(x)\leq c\}$.
Recall that $\cl V\subset U$ in view of \cref{L2.1}\,\ref{L2.1.5}.
By the assumptions, $\{x_{1k}\}$ has accumulation points, and they all belong to $U$.
For any accumulation point $\hat x\in U$ and the corresponding subsequence
(without relabeling) $\{x_{1k}\}$ with $x_{1k}\to\hat x$, we also have
$x_{2k}\to\hat x$.
Thanks to the lower semicontinuity of $\varphi_1$ and $\varphi_2$, it holds
$(\varphi_1+\varphi_2)(\hat x)\le \liminf_{k\to+\infty} (\varphi_1(x_{1k})+\varphi_1(x_{2k}))\le\al$, and consequently, condition \eqref{ULCb} is satisfied
with $x_k:=\hat x$ for all $k\in\N$.
In all cases, the conclusions follow from \cref{P4.3}.
\end{proof}

\begin{remark}
\begin{enumerate}
\item
\cref{P4.6}\,\ref{P4.6.1} strengthens \cite[Proposition~2.7.3]{BorweinZhu96}, \cite[Proposition~2.1\,(a) and (c)]{Las01} and \cite[Lemma~1.123\,(c)]{Pen13}.

\item
If $U$ is compact, then all the compactness assumptions as well as assumption
$\inf_U\varphi_2>-\infty$ in \cref{P4.6} are satisfied automatically.

\item
As illustrated by \cref{E3.4} (see also \cref{R4.2}\,\ref{R4.2.7}), condition \eqref{P4.6-4} in part \ref{P4.6.2} of \cref{P4.6} is essential.

\item
Using slightly weaker versions of the `firm' properties, corresponding to replacing
${\Theta}_{U}^\circ(\varphi_1,\varphi_2)$ and
${\Theta}_{U}^\dag(\varphi_1,\varphi_2)$ in parts \ref{D4.1.3} and \ref{D4.1.4} of \cref{D4.1}
by the \RHS s of the respective inequalities in \cref{P3.2}\,\ref{P3.2.3}
(see \cref{R4.1}\,\ref{R4.1.5}), one can drop assumptions \eqref{P4.6-3} and \eqref{P4.6-4} in parts \ref{P4.6.1} and \ref{P4.6.2} of \cref{P4.6},
thus, strengthening these assertions.
\end{enumerate}
\end{remark}

When it comes down to minimization
problems in infinite-dimensional spaces,
weakly sequentially lower semicontinuous functions are of special interest.

\begin{proposition}
\label{P4.7}
Let $X$ be a normed space,
$\varphi_1$ and $\varphi_2$ be weakly sequentially \lsc\ on $U$ and $\inf_U\varphi_2>-\infty$.
The pair $(\varphi_1,\varphi_2)$ is \begin{enumerate}
\item
{uniformly lower semicontinuous}
on $U$ if $\{x\in U\,|\,\varphi_1(x)\leq c\}$ is weakly sequentially compact for each $c\in\R$;

\item
\label{P4.7.2}
{quasiuniformly lower semicontinuous}
on $U$ if $\cl^wV\subset U$ and $\{x\in\cl^w V\,|\,\varphi_1(x)\leq c\}$ is weakly sequentially compact for each $V\in EI(U)$ and $c\in\R$,
where $\cl^wV$ stands for the weak closure of $V$, i.e., the closure of $V$ with respect to the weak topology in $X$.
\end{enumerate}
\end{proposition}

\begin{proof}[Sketch of the proof]
The proof basically repeats that of the corresponding parts of \cref{P4.6},
replacing strong convergence with the weak one.
If $\{x_{1k}\}$ and $\{x_{2k}\}$ are the sequences constructed in the above proof, then,
by the weak sequential compactness assumptions, $\{x_{1k}\}$ has weak accumulation points.
Thanks to the weak sequential lower semicontinuity of $\varphi_1$ and $\varphi_2$, it holds
$(\varphi_1+\varphi_2)(\hat x)\le \liminf_{k\to+\infty} (\varphi_1(x_{1k})+\varphi_1(x_{2k}))$, and the conclusions follow from \cref{P4.3}.
\end{proof}

Note that weak convergence does not allow us to establish similar sufficient conditions for the `firm' properties.

\begin{remark}
The assertion in part \ref{P4.7.2} of \cref{P4.7} remains true if $\cl^w V$ is replaced by the
weak \emph{sequential} closure of $V$. However, the weak sequential closure of a set does not need to be
weakly sequentially closed in general, see e.g.\ \cite{MehlitzWachsmuth2019} for a study addressing so-called decomposable sets
in Lebesgue spaces, so the associated statement has to be used with care.
\end{remark}

\begin{corollary}
\label{C4.3}
Let $X$ be a reflexive Banach space, $U$ be convex and bounded and $\varphi_1$ and $\varphi_2$ be weakly sequentially \lsc\ on~$U$.
Then the pair $(\varphi_1,\varphi_2)$ is quasiuniformly lower semicontinuous on $U$.
\end{corollary}

\begin{proof}
Let $V\in EI(U)$.
By \cref{L2.1}\,\ref{L2.1.10}, $\co V\in EI(U)$, and, by \cref{L2.1}\,\ref{L2.1.5}, $\cl\co V\in EI(U)$.
Hence, $\cl^wV\subset\cl^w\co V=\cl\co V\subset U$.
Since $\cl^w V$ is, particularly, weakly sequentially closed and bounded while the sublevel sets
of weakly sequentially lower semicontinuous functions are weakly sequentially closed,
the set $\{x\in\cl^w V\,|\,\varphi_1(x)\leq c\}$ is weakly sequentially compact for each $c\in\R$ as $X$ is reflexive.
The assertion follows from \cref{P4.7}\,\ref{P4.7.2}.
\end{proof}

\section{Relative uniform and quasiuniform lower semicontinuity}
\label{S5}

In this section we investigate the situation where at least one of the involved functions is the indicator function of a set.

We start with the case $\varphi_1:=\varphi$ and $\varphi_2:=i_\Omega$ for some function
$\varphi\colon X\to\R_\infty$ on a metric space $X$ and subset $\Omega\subset X$.
Our basic assumption \eqref{nonempty} in this setting becomes
\begin{gather}
\label{nonempty2}
\dom\varphi\cap\Omega\cap U\ne\es,
\end{gather}
where $U$ is another subset of $X$.
From \eqref{La0}, \eqref{Th}, \eqref{phihat},
\cref{P3.2}\,\ref{P3.2.4} and \cref{R3.2}, we obtain:
\begin{subequations}
\label{LaMOm}
\begin{align}
\label{LaOm0}
\Lambda_U^\circ(\varphi,i_\Omega)=& \liminf_{\substack{\dist(x,\Omega\cap U)\to0,\, x\in U}} \varphi(x),
\\
\label{MOm0}
{\Lambda}_U^\dag(\varphi,i_\Omega)=& \inf_{V\in EI(U)}\; 
\Lambda_V^\circ(\varphi,i_\Omega),
\\
\label{ThOm0}
{\Theta}_{U}^\circ(\varphi,i_\Omega)=&
\limsup_{\substack{\dist(x,\Omega\cap U)\to0\\x\in\dom\varphi\cap U}}\; \inf\limits_{\substack{u\in\Omega\cap U}}
\max\{d(u,x),\varphi(u)-\varphi(x)\},
\\
\label{MOmhat}
{\Theta}_{U}^\dag(\varphi,i_\Omega)=&
\sup\limits_{\substack{V\in EI(U)}}\;
\limsup\limits_{\substack{\dist(x,\Omega)\to0\\ x\in\dom\varphi\cap V}}\;
\inf\limits_{\substack{u\in\Omega\cap U}}
\max\{d(u,x),\varphi(u)-\varphi(x)\}.
\end{align}
\end{subequations}
When $U=X$, quantities \eqref{LaOm0} and \eqref{MOm0} coincide (see \cref{P3.1}\,\ref{P3.1.2})
and reduce to the \emph{uniform infimum} \cite[page~1034]{Las01} (\emph{decoupled infimum} \cite[Definition~3.2.1]{BorZhu05},
\emph{stabilized infimum} \cite[Definition~1.122]{Pen13})
of $\varphi$ on~$\Omega$.

The uniform and quasiuniform lower semicontinuity properties in \cref{D4.1}
reduce to the corresponding ones in the next definition employing \eqref{LaMOm}.

\begin{definition}
\label{D5.1}
The function $\varphi$ is
\begin{enumerate}
\item
\label{D5.1.1}
\emph{uniformly lower semicontinuous relative to} $\Omega$
on $U$ if $\inf_{\Omega\cap U} \varphi\le{\Lambda}_U^\circ(\varphi,i_\Omega)$;
\item
\label{D5.1.2}
\emph{quasiuniformly lower semicontinuous relative to} $\Omega$
on $U$ if $\inf_{\Omega\cap U} \varphi\le{\Lambda}_U^\dag(\varphi,i_\Omega)$;
\item
\label{D5.1.3}
\emph{firmly uniformly lower semicontinuous relative to} $\Omega$
on $U$ if ${\Theta}_{U}^\circ(\varphi,i_\Omega)=0$;
\item
\label{D5.1.4}
\emph{firmly quasiuniformly lower semicontinuous relative to} $\Omega$
on $U$ if ${\Theta}_{U}^\dag(\varphi,i_\Omega)=0$;
\item
uniformly/quasiuniformly/firmly uniformly/firmly quasiuniformly lower semicontinuous relative to $\Omega$
near a point $\bx\in\dom\varphi\cap\Omega$
if it is uniformly/quasiuniformly/firmly uniformly/firmly quasiuniformly lower semicontinuous
relative to $\Omega$ on $\overline B_\de(\bx)$
for all sufficiently small $\de>0$.
\end{enumerate}
\end{definition}

\begin{remark}
\begin{enumerate}
\item
The quasiuniform relative lower semicontinuity properties in \cref{D5.1} recapture the corresponding relative lower semicontinuity properties in \cite[Definition~3.3]{KruMeh22}.

\item
In view of \cref{R4.1}\,\ref{R4.1.2}
the properties in \cref{D5.1}\,\ref{D5.1.1}--\ref{D5.1.4} are related by \eqref{R4.1-2}.
\end{enumerate}
\end{remark}

The characterizations of the relative lower semicontinuity properties in the next two propositions
are consequences of representations \eqref{LaMOm} and corresponding assertions in \cref{P3.1,P3.2,P3.3}.
They can also be derived from \cref{P4.1}.

\begin{proposition}
\label{P5.1}
The function $\varphi$ is
\begin{enumerate}
\item
\label{P5.1.1}
uniformly lower semicontinuous relative to $\Omega$ on $U$ if and only if, for any $\eps>0$, there exists an $\eta>0$ such that, for any
$x\in U$
with $\dist(x,\Omega\cap U)<\eta$, there is a $u\in\Omega\cap U$ such that
$\varphi(u)<\varphi(x)+\eps$;

\item
quasiuniformly lower semicontinuous relative to $\Omega$ on $U$ if and only if, for any $V\in EI(U)$ and $\eps>0$, there exists an $\eta>0$ such that, for any
$x\in V$
with $\dist(x,\Omega)<\eta$, there is a $u\in\Omega\cap U$ such that
$\varphi(u)<\varphi(x)+\eps$;

\item
\label{P5.1.3}
{firmly uniformly lower semicontinuous relative to} $\Omega$ on $U$ if and only if, for any $\eps>0$, there exists an $\eta>0$ such that, for any
$x\in\dom\varphi\cap U$
with $\dist(x,\Omega\cap U)<\eta$, there is a $u\in\Omega\cap U\cap B_\eps(x)$ such that
$\varphi(u)<\varphi(x)+\eps$;

\item
\label{P5.1.4}
{firmly quasiuniformly lower semicontinuous relative to} $\Omega$ on $U$ if and only if, for any $V\in EI(U)$ and $\eps>0$, there exists an $\eta>0$ such that, for any
$x\in\dom\varphi\cap V$
with $\dist(x,\Omega)<\eta$, there is a $u\in\Omega\cap B_\eps(x)$ such that
$\varphi(u)<\varphi(x)+\eps$.
\end{enumerate}
\end{proposition}

\begin{remark}
\cref{P5.1}\,\ref{P5.1.1} strengthens \cite[Lemma~1.123\,(a)]{Pen13}:
if $\varphi$ is uniformly lower semicontinuous \emph{around} $\Omega$ \cite[page 88]{Pen13}, then it is uniformly lower semicontinuous relative to $\Omega$ on~$X$.
\end{remark}

The next proposition gives a simplified characterization of
firm quasiuniform relative lower semicontinuity near a point.
It is a consequence of \cref{prop:firm_quasiuniform_lsc_near_point}.

\begin{proposition}
\label{P5.2}
	The function $\varphi$ is
	firmly quasiuniformly lower semicontinuous
	near $\bx\in\dom\varphi\cap\Omega$ if and only if $\Theta^\dag_{B_\de(\bar x)}(\varphi,i_\Omega)=0$ for some $\de>0$.
\end{proposition}

\cref{E3.1} illustrates the firm uniform relative lower semicontinuity property (see \cref{R4.2}\,\ref{R4.2.7}).
The following infinite-dimensional example will be important later on when we
discuss applications of our findings.
\sloppy

\begin{example}
\label{ex:lower_semicontinuity_relative_to_a_set}
	Let $D\subset\R^d$ be a Lebesgue-measurable set with positive and finite Lebesgue measure $\blambda(D)$.
	We equip $D$ with the $\sigma$-algebra of all Lebesgue-measurable subsets of $D$ as well
	as (the associated restriction of) the Lebesgue measure $\blambda$,
and consider the Lebesgue space $L^2(D)$ of all
	(equivalence classes of) measurable, square integrable functions equipped with the
	usual norm.
In what follows, we suppress
\emph{Lebesgue} for brevity.
	
Define a function $\varphi\colon L^2(D)\to\R$ by means of
\begin{equation}\label{eq:sparsity_promoting_function}
		\forall x\in L^2(D)\colon\quad
		\varphi(x):=\blambda(\{x\neq 0\}).
	\end{equation}
We use the notation
$\{x\neq 0\}:=\{\omega\in D\,|\,x(\omega)\neq 0\}$ for brevity.
Furthermore, the sets $\{x=0\}$, $\{x<0\}$, $\{x>0\}$ and analogous sets
with non-vanishing right-hand side or bilateral bounds
are defined similarly.
We note that, by definition of $L^2(D)$, these sets are well defined
up to subsets of measure zero. Particularly, $\varphi$ from
\eqref{eq:sparsity_promoting_function} is well defined.
By means of Fatou's lemma one can easily check that $\varphi$
	is lower semicontinuous, see \cite[Lemma~2.2]{MehlitzWachsmuth22}.
	
	For fixed functions $x_a,x_b\in L^2(D)$ satisfying $x_a(\omega)\leq 0\leq x_b(\omega)$
	for almost all $\omega\in D$, we define the
\emph{box-constraint} set
	$\Omega\subset L^2(D)$ by means of
	\begin{equation}\label{eq:box_constrained_set}
		\Omega:=\{x\in L^2(\Omega)\,|\,x_a\leq x\leq x_b\text{ a.e.\ on }D\},
	\end{equation}
and note that $\Omega$ is nonempty, closed and convex.
For an $x\in L^2(D)$, we define
	\begin{align*}
		\forall \omega\in D\colon\quad
		u_x(\omega)
		&:=
		\begin{cases}
			x_a(\omega)	&	\omega\in \{x<x_a\},\\
			x_b(\omega)	&	\omega\in \{x>x_b\},\\
			x(\omega)	&	\omega\in \{x_a\leq x\leq x_b\}.
		\end{cases}
	\end{align*}
A simple calculation shows that
	$u_x\in L^2(D)$ is the uniquely determined projection of $x$ onto $\Omega$,
and consequently, $\dist(x,\Omega)=\norm{x-u_x}$.
Furthermore, by construction, we have
$\varphi(u_x)\leq\varphi(x)$.
Given any $\bar x\in\Omega$, $\de>0$, $\eps>0$ and $\eta\in(0,\eps)$,
conditions $x\in  B_\de(\bar x)$ and $\dist(x,\Omega)<\eta$
yield $u_x\in\Omega\cap B_\eps(x)$ and
\begin{align*}
	\norm{u_x-\bar x}^2
	&=
	\int_{\{x<x_a\}}(x_a-\bar x)^2\,\mathrm d\omega
	+
	\int_{\{x_a\leq x\leq x_b\}}(x-\bar x)^2\,\mathrm d\omega
	+
	\int_{\{x<x_b\}}(x_b-\bar x)^2\,\mathrm d\omega
	\\
	&
	\leq
	\int_{\{x<x_a\}}(x-\bar x)^2\,\mathrm d\omega
	+
	\int_{\{x_a\leq x\leq x_b\}}(x-\bar x)^2\,\mathrm d\omega
	+
	\int_{\{x<x_b\}}(x-\bar x)^2\,\mathrm d\omega
	\\
	&=
	\norm{x-\bar x}^2
	\le
	\delta^2,
\end{align*}
i.e., $u_x\in B_\de(\bar x)$.
By \cref{R4.2}\,\ref{R4.2.5} and \cref{P5.1}\,\ref{P5.1.3},
$\varphi$ is firmly uniformly lower semicontinuous relative to $\Omega$ near any point in $\Omega$.
	
	Note that the function $\varphi$ is discontinuous and not weakly sequentially lower semicontinuous.
	In fact, $\varphi$ is \emph{nowhere} Lipschitz continuous,
	see \cite[Corollary~3.9]{MehlitzWachsmuth22}.
	Clearly, $i_\Omega$ is discontinuous.
	Thus, we have constructed a uniformly lower semicontinuous pair of
	functions, both non-Lipschitz, while one of them is not weakly sequentially
	lower semicontinuous.
	
\end{example}

Next, we discuss sufficient conditions for firm uniform and quasiuniform lower semicontinuity of a function relative to a set.
For sufficient conditions for (not firm) quasiuniform lower semicontinuity of a function relative
to a set, we refer the interested reader to \cite[Section~3.3]{KruMeh22}.
The next two statements are direct consequences of \cref{P4.5,P4.6,P5.1}.

\begin{proposition}
\label{P5.3}
If $\varphi$ is uniformly continuous on $U$, then
it is firmly uniformly lower semicontinuous relative to $\Omega$ on $U$.
\end{proposition}

\begin{proposition}
\label{P5.4}
Let $\varphi$ be \lsc\ on $U$ and $\Omega$ be closed.
The function $\varphi$ is
\begin{enumerate}
\item
firmly uniformly lower semicontinuous relative to $\Omega$ on $U$ if $\{x\in U\,|\,\varphi(x)\leq c\}$ is compact for each $c\in\R$, and
\[
	\limsup\limits_{\substack{\dist(x,\Omega\cap U)\to0,\, x\in\dom\varphi\cap U}}\; \varphi(x)<+\infty;
\]
\item
firmly quasiuniformly lower semicontinuous relative to $\Omega$ on $U$ if
$\{x\in\cl V\,|\,\varphi(x)\leq c\}$ is compact for each $V\in EI(U)$ and $c\in\R$, and
\[
\sup\limits_{\substack{V\in EI(U)}}\; \limsup\limits_{\substack{\dist(x,\Omega)\to0,\, x\in\dom\varphi\cap V}}\; \varphi(x)<+\infty.
\]
\end{enumerate}
\end{proposition}

The next example illustrates the difference between uniform and quasiuniform relative lower semicontinuity.

\begin{example}
Let closed convex sets $\Omega,U\subset\R^2$ and a convex function $\varphi\colon\R^2\to\R_\infty$ be given by
\begin{gather*}
\Omega:=\{(x,y)\mid y\le0\},
\quad
U:=\{(x,y)\mid y\ge x^2\},
\\
\forall (x,y)\in\R^2\colon\quad
\varphi(x,y):=
\begin{cases}
-1&\text{if }y>0,
\\
0&\text{if }y=0,
\\
+\infty&\text{otherwise}.
\end{cases}
\end{gather*}
Then $\dom\varphi\cap U=U$, $\Omega\cap U=\{(0,0)\}$,
\begin{gather*}
\inf_{\Omega\cap U}\varphi=\varphi(0,0)=0\AND
\liminf_{\substack{\dist((x,y),\Omega)\to0,\; (x,y)\in U}} \varphi(x,y)=-1.
\end{gather*}
By \cref{P5.1}\,\ref{P5.1.1} $\varphi$ is not uniformly lower semicontinuous relative to $\Omega$ on $U$.
At the same time $\dist(V,\Omega)>0$ for any $V\in EI(U)$.
Hence, if $\eta\in(0,\dist(V,\Omega))$, then $\dist(x,\Omega)>\eta$ for all $x\in V$,
and the conditions in \cref{P5.1}\,\ref{P5.1.4} are trivially satisfied, i.e., $\varphi$ is firmly quasiuniformly lower semicontinuous relative to $\Omega$ on $U$.
Observe that $\varphi$ is not lower semicontinuous at~$(0,0)$.
\end{example}


The case of two indicator functions of some subsets $\Omega_1,\Omega_2\subset X$ can be considered as a particular case of the uniform/quasiuniform lower semicontinuity properties in \cref{D4.1} or relative lower semicontinuity properties in \cref{D5.1}.
The corresponding properties are rather weak and are satisfied almost automatically in most natural situations.

Let $\Omega_1,\Omega_2,U\subset X$ and  $\Omega_1\cap\Omega_2\cap U\ne\es$.
First, observe that $\inf_U(i_{\Omega_1}+i_{\Omega_2})= {\Lambda}_U^\circ(i_{\Omega_1},i_{\Omega_2})=0$.
Hence, $(i_{\Omega_1},i_{\Omega_2})$ is automatically uniformly (and quasiuniformly) lower semicontinuous on $U$.
Using \eqref{ThOm0}, \eqref{MOmhat} and \cref{P5.2} as well as parts \ref{P3.3.1} and \ref{P3.3.4} of \cref{P3.3}, we can formulate characterizations of firm uniform and quasiuniform lower semicontinuity.

\begin{proposition}
\label{P5.5}
The pair $(i_{\Omega_1},i_{\Omega_2})$ is
\begin{enumerate}
\item
\label{P5.5.1}
firmly uniformly lower semicontinuous on $U$ if and only if
\begin{align}
\label{P5.5-1}
\limsup_{\substack{\dist(x,\Omega_2\cap U)\to0,\, x\in\Omega_1\cap U}}\; \dist(x,\Omega_1\cap\Omega_2\cap U)=0;
\end{align}
\item
\label{P5.5.2}
firmly quasiuniformly lower semicontinuous on $U$ if and only if
\begin{gather}
\label{P5.5-2}
\sup\limits_{V\in EI(U)}\;
\limsup\limits_{\substack{\dist(x,\Omega_2)\to 0,\; x\in\Omega_1\cap V}}\;
\dist(x,\Omega_1\cap\Omega_2)=0;
\end{gather}
\item
firmly quasiuniformly lower semicontinuous near a point $\bx\in\Omega_1\cap\Omega_2$ if and only if for some $\de>0$ it holds:
\begin{gather*}
\limsup\limits_{\substack{\dist(x,\Omega_2\cap B_\de(\bx))\to 0,\; x\in\Omega_1\cap B_\de(\bx)}}\;
\dist(x,\Omega_1\cap\Omega_2)=0.
\end{gather*}
\end{enumerate}
\end{proposition}

The next proposition is a consequence of \cref{P5.4}.
The statement
and its corollary show that the situations when a pair of indicator functions is not firmly uniformly or firmly quasiuniformly lower semicontinuous are rare.

\begin{proposition}
\label{P5.6}
Let $\Omega_1$ and $\Omega_2$ be closed.
The pair $(i_{\Omega_1},i_{\Omega_2})$ is
\begin{enumerate}
\item
firmly uniformly lower semicontinuous relative to $\Omega$ on $U$ if $\Omega_1\cap U$ is compact;
\item
firmly quasiuniformly lower semicontinuous on $U$ if the sets $\Omega_1\cap\cl V$ are compact for all $V\in EI(U)$.
\end{enumerate}
\end{proposition}

\begin{corollary}
\label{C5.4}
Let $X$ be a finite dimensional Banach space,
$\Omega_1,\Omega_2$ be closed, and $U$ be bounded.
Then $(i_{\Omega_1},i_{\Omega_2})$ is firmly quasiuniformly lower semicontinuous on $U$.
\end{corollary}

The next statement gives alternative characterizations of the firm uniform and quasiuniform lower semicontinuity of a pair of indicator functions.

\begin{proposition}
\label{P5.7}
The pair $(i_{\Omega_1},i_{\Omega_2})$ is
\begin{enumerate}
\item
firmly uniformly lower semicontinuous on $U$ if and only if
\begin{gather}
\label{P5.7-2}
\limsup_{\substack{\dist(x,\Omega_1\cap U)\to0,\; \dist(x,\Omega_2\cap U)\to0,\; x\in U}} \dist(x,\Omega_1\cap\Omega_2\cap U)=0;
\end{gather}
\item
firmly quasiuniformly lower semicontinuous on $U$ if and only if
\begin{gather}
\label{P5.7-1}
\sup_{V\in EI(U)}\; \limsup_{\substack{\dist(x,\Omega_1)\to0,\; \dist(x,\Omega_2)\to0,\; x\in V}} \dist(x,\Omega_1\cap\Omega_2)=0;
\end{gather}
\item
\label{P5.7.3}
firmly quasiuniformly lower semicontinuous near a point $\bx\in\Omega_1\cap\Omega_2$ if and only if for some $\de>0$ it holds:
\begin{gather}
\label{P5.7-3}
\limsup\limits_{\substack{\dist(x,\Omega_1)\to 0,\; \dist(x,\Omega_2)\to 0,\; x\in B_\de(\bx)}}\;
\dist(x,\Omega_1\cap\Omega_2)=0.
\end{gather}
\end{enumerate}
\end{proposition}

\begin{proof}
We prove the second assertion.
The proofs of the first and the third ones follow the same pattern with some obvious simplifications.
Observe that \eqref{P5.7-1} trivially implies \eqref{P5.5-2}.
Thanks to \cref{P5.5}\,\ref{P5.5.2}, it suffices to show the opposite implication.
Let condition \eqref{P5.5-2} be satisfied, and let
$V\in EI(U)$,
$\{x_k\}\subset V$,
$\dist(x_k,\Omega_1)\to0$ and
$\dist(x_k,\Omega_2)\to0$ as $k\to+\infty$.
Then for each $k\in\N$, there exist points $x_{1k}\in\Omega_1$ and $x_{2k}\in\Omega_2$ such that $d(x_{1k},x_k)\to0$ and $d(x_{2k},x_k)\to0$; hence, $\dist(x_{1k},\Omega_2)\le d(x_{1k},x_{2k})\le d(x_{1k},x_{k})+d(x_{2k},x_{k})\to0$.
By \cref{L2.1}\,\ref{L2.1.6}, there exists a subset $V'\in EI(U)$ such that $V\in EI(V')$.
Then $x_{1k}\in V'$ for all sufficiently large $k\in\N$.
By \eqref{P5.5-2},
$\dist(x_{1k},\Omega_1\cap\Omega_2)\to0$, and consequently,
$\dist(x_{k},\Omega_1\cap\Omega_2)\to0$.
Thus, condition \eqref{P5.7-1} holds true.
\end{proof}

\begin{remark}
In view of \cref{P5.7}\,\ref{P5.7.3},
the firm quasiuniform lower semicontinuity of a pair of indicator functions near a point in the intersection of the sets is implied,
for instance, by the well known and widely used \emph{subtransversality} property
(also known as \emph{linear regularity, metric regularity, linear coherence} and \emph{metric inequality}),
and as a consequence, also by the stronger \emph{transversality} property (also known under various names);
see, e.g., \cite{Iof17,KruLukTha18,BuiCuoKru20}.
Recall that the sets $\Omega_1$ and $\Omega_2$ are \emph{subtransversal} at $\bx\in\Omega_1\cap\Omega_2$
if there exist numbers $\al>0$ and $\de>0$ such that
\begin{gather}
\label{R5.3-1}
\forall x\in B_\de(\bar x)\colon\quad
\dist(x,\Omega_1\cap\Omega_2)\le \al\max\{\dist(x,\Omega_1),\dist(x,\Omega_2)\}.
\end{gather}
Nonlocal versions of this property, i.e., with some subset $U\subset X$ (e.g., $U=X$) in place of $B_\de(\bx)$ are also in use.
Condition \eqref{R5.3-1} describes so called \emph{linear subtransversality}.
More subtle nonlinear, in particular, H\"older subtransversality (see e.g.\ \cite{CuoKru20.2}) is still sufficient for the property~\eqref{P5.7-3}.
\end{remark}

The following example, which is inspired by \cref{E3.1}, shows that the firm uniform lower semicontinuity of a pair of indicator functions
can be strictly weaker than (linear) subtransversality of the involved sets.

\begin{example}
\label{E5.3}
Let closed convex sets $\Omega_1,\Omega_2\subset\R^2$ be given by
\begin{align*}
\Omega_1:=\{(x,y)\in\R^2\mid y\ge x^2\},
\quad
\Omega_2:=\{(x,y)\in\R^2\mid y\le0\}.
\end{align*}
Then $\Omega_1\cap\Omega_2=\{(0,0)\}$.
If $(x,y)\in\Omega_1$ and
$\dist((x,y),\Omega_2)\to0$, then $(x,y)\to(0,0)$.
Hence, given any subset $U\subset\R^2$ containing $(0,0)$, we have $\dist((x,y),\Omega_1\cap\Omega_2\cap U)\to0$,
i.e., condition \eqref{P5.5-1} is satisfied, and, by \cref{P5.5}\,\ref{P5.5.1}, the pair $(i_{\Omega_1},i_{\Omega_2})$ is
firmly uniformly lower semicontinuous on~$U$.
At the same time, considering the points $(x_k,y_k):=(1/k,1/k^2)$ as $k\to+\infty$,
one can easily check that condition \eqref{R5.3-1} fails for any $\al>0$ and $\de>0$, i.e.,
$\Omega_1$ and $\Omega_2$ are not subtransversal at $(0,0)$.
\end{example}

\section{Optimality conditions}
\label{S6}

We consider here the problem of minimizing the sum of two functions $\varphi_1,\varphi_2\colon X\to\R_\infty$ on a metric space $X$.
When discussing dual optimality conditions, $X$ will be assumed Banach or, more specifically, Asplund.
This model is quite general (see a discussion in \cite{Las01}).
It may represent so-called \emph{composite optimization} problems, where typically the
smoothness properties of $\varphi_1$ and $\varphi_2$ are rather different.
If one of the functions is an indicator function, the model covers
nonsmooth constrained optimization problems.
As in the previous sections, we are going to exploit the decoupling approach, allowing $\varphi_1$ and $\varphi_2$ to take different inputs.

We mostly discuss local minimality/stationarity properties of $\varphi_1+\varphi_2$ around a given point $\bx\in\dom\varphi_1\cap\dom\varphi_2$.
Recall that $\bx$ is called a \emph{local uniform minimum} \cite{Las01} of $\varphi_1+\varphi_2$ if it satisfies \eqref{UM}.
This notion is stronger than the conventional local minimum.
Together with the related definitions of \emph{uniform infimum} \eqref{La} and \eqref{La0} and \emph{uniform lower semicontinuity} \eqref{La0qc}
they form the foundations of the decoupling approach; see \cite{BorweinZhu96,BorweinIoffe1996,Las01,BorZhu05}.
In what follows, we examine weaker local \emph{quasiuniform} minimality and stationarity concepts
which are based on the decoupling quantity \eqref{M}, and utilise the properties discussed in \cref{S3} as well as the \emph{quasiuniform lower semicontinuity} from \cref{S4}.

\begin{definition}
\label{D6.1}
\begin{enumerate}
\item
The point $\bx$ is a \emph{local quasiuniform minimum} of $\varphi_1+\varphi_2$
if condition \eqref{QUM} is satisfied.
If the latter condition is satisfied with $\de=+\infty$, then
$\bx$ is referred to as a \emph{quasiuniform minimum} of $\varphi_1+\varphi_2$.
\item
\label{D6.1.3}
Given $\eps>0$ and $\de>0$,
$\bx$ is a \emph{quasiuniform $\eps$-minimum} of $\varphi_1+\varphi_2$ on $B_\de(\bx)$ if
\begin{align}
\label{D6.1-3}
(\varphi_1+\varphi_2)(\bx)< {\Lambda}_{B_\de(\bx)}^\dag(\varphi_1,\varphi_2)+\eps.
\end{align}

\item
\label{D6.1.2}
The point $\bx$ is a \emph{quasiuniform stationary} point of $\varphi_1+\varphi_2$ if
for any $\varepsilon>0$, there exists a $\de_\varepsilon>0$ such that,
for any $\de\in(0,\de_\varepsilon)$, $\bar x$ is a quasiuniform $\varepsilon\de$-mini\-mum of $\varphi_1+\varphi_2$ on $B_\de(\bar x)$.
\end{enumerate}
\end{definition}

\begin{proposition}
Every local uniform minimum of $\varphi_1+\varphi_2$ is a local quasiuniform minimum of $\varphi_1+\varphi_2$.
\end{proposition}

\begin{proof}
As observed in \cref{sec:introduction}, expression $\Lambda_{B_\de(\bx)}(\varphi_1,\varphi_2)$ in definition \eqref{UM} of local uniform minimum can be replaced with
$\Lambda_{B_\de(\bx)}^\circ(\varphi_1,\varphi_2)$ (this is an immediate consequence of assertions \ref{P3.1.1}, \ref{P3.1.3} and \ref{P3.1.6} of \cref{P3.1}).
Thanks to this observation, the assertion is a consequence of \cref{P3.1}\,\ref{P3.1.4}.
\end{proof}

\begin{proposition}
\label{prop:quasiuniform_min_vs_stat}
Every local quasiuniform minimum of $\varphi_1+\varphi_2$ is a quasiuniform stationary point of $\varphi_1+\varphi_2$.
\end{proposition}

\begin{proof}
Let $\bx$ be a local quasiuniform minimum of $\varphi_1+\varphi_2$.
Then, by \eqref{QUM}, there is a $\bar\de>0$ such that $(\varphi_1+\varphi_2)(\bx)= {\Lambda}_{B_{\bar\de}(\bx)}^\dag(\varphi_1,\varphi_2)$.
Given any $\eps>0$ and $\de\in(0,\bar\de)$, in view of \cref{P3.1}\,\ref{P3.1.3}, we have
$$(\varphi_1+\varphi_2)(\bx)\le {\Lambda}_{B_\de(\bx)}^\dag(\varphi_1,\varphi_2)< {\Lambda}_{B_\de(\bx)}^\dag(\varphi_1,\varphi_2)+\eps\de.$$
Hence, $\bx$ is a quasiuniform stationary point.
\end{proof}

The properties in \cref{D6.1} imply the corresponding
conventional local minimality/stationarity properties of $\varphi_1+\varphi_2$.
As the following proposition reveals, they become equivalent when the pair $(\varphi_1,\varphi_2)$ is quasiuniformly lower semicontinuous
(in the sense of \cref{D4.1}\,\ref{D4.1.2}) on an appropriate \nbh\ of $\bx$.

\begin{proposition}
\label{P6.3}
\begin{enumerate}
\item
\label{P6.3.1}
If $\bar x$ is a local quasiuniform minimum of $\varphi_1+\varphi_2$, then it is a local minimum of $\varphi_1+\varphi_2$.
The two properties are equivalent
provided that $(\varphi_1,\varphi_2)$ is quasiuniformly lower semicontinuous near $\bx$.
\item
\label{P6.3.2}
Let $\varepsilon>0$ and $\de>0$.
If $\bar x$ is a quasiuniform $\varepsilon$-minimum of $\varphi_1+\varphi_2$
on $B_\de(\bar x)$, then it is an $\varepsilon$-minimum of $\varphi_1+\varphi_2$ on $B_\de(\bar x)$.
The two properties are equivalent
provided that $(\varphi_1,\varphi_2)$ is quasiuniformly lower semicontinuous on $B_\de(\bx)$.
\item
\label{P6.3.3}
If $\bar x$ is a quasiuniform stationary point of $\varphi_1+\varphi_2$, then $\bar x$ is stationary for $\varphi_1+\varphi_2$.
The two properties are equivalent
provided that $(\varphi_1,\varphi_2)$ is quasiuniformly lower semicontinuous near $\bx$.
\end{enumerate}
\end{proposition}

\begin{proof}
The `if' parts of assertions \ref{P6.3.1} and \ref{P6.3.2} are direct consequences of \cref{P3.1}\,\ref{P3.1.4}.
The `if' part of assertion \ref{P6.3.3} is a consequence of \cref{P3.1}\,\ref{P3.1.4} and \Cref{lem:characterization_stationarity}.

We now prove the converse implication in assertion \ref{P6.3.2}.
Let $(\varphi_1,\varphi_2)$ be quasiuniformly lower semicontinuous on $B_\de(\bx)$ and
$\bx$ be an $\eps$-minimum of $\varphi_1+\varphi_2$ on ${B}_{\de}(\bx)$.
Then it is an $\eps'$-minimum of $\varphi_1+\varphi_2$ on ${B}_{\de}(\bx)$ for some $\eps'\in(0,\eps)$.
Choose a number $\xi\in(0,\eps-\eps')$ and a subset $V\in EI({B}_{\de}(\bx))$.
By \cref{P4.1}\,\ref{P4.1.2}, there exists an $\eta>0$ such that, for any
$x_1\in V$ and $x_2\in X$
with $d(x_1,x_2)<\eta$, there is an $x\in B_\de(\bx)$ satisfying
$(\varphi_1+\varphi_2)(x)< \varphi_1(x_{1})+\varphi_2(x_{2})+\xi$,
and consequently,
\begin{align*}
(\varphi_1+\varphi_2)(\bar x)<
(\varphi_1+\varphi_2)(x)+\eps'<
\varphi_1(x_{1})+\varphi_2(x_{2})+\eps'+\xi.
\end{align*}
Hence,
\begin{align*}
(\varphi_1+\varphi_2)(\bar x)\le
\inf_{V\in EI({B}_{\de}(\bx))}\; \liminf_{d(x_1,x_2)\to0,\, x_1\in V}
(\varphi_1(x_{1})+\varphi_2(x_{2}))+\eps'+\xi.
\end{align*}
Since $\eps'+\xi<\eps$, and in view of \cref{P3.1}\,\ref{P3.1.7} and \cref{D6.1}\,\ref{D6.1.3}, $\bx$ is a quasiuniform $\eps$-minimum of $\varphi_1+\varphi_2$ on $B_\de(\bx)$.

For the remainder of the proof, let $(\varphi_1,\varphi_2)$ be quasiuniformly lower semicontinuous near $\bx$.
Suppose that $\bar x$ is a local minimum of $\varphi_1+\varphi_2$, and
let a $\de>0$ be such that $(\varphi_1,\varphi_2)$ is quasiuniformly lower semicontinuous on $B_\de(\bx)$ and
$\bx$ is a minimum of $\varphi_1+\varphi_2$ on $B_\de(\bx)$.
Then $\bx$ is also an $\eps$-minimum of $\varphi_1+\varphi_2$ on $B_\de(\bx)$ for any $\eps>0$.
As shown above, $\bx$ is a quasiuniform $\eps$-minimum of $\varphi_1+\varphi_2$ on $B_\de(\bx)$, i.e., condition \eqref{D6.1-3} is satisfied.
Letting $\eps\downarrow0$, we arrive at
$(\varphi_1+\varphi_2)(\bx)\le {\Lambda}_{B_\de(\bx)}^\dag (\varphi_1,\varphi_2)$.
Since, due to \cref{P3.1}\,\ref{P3.1.4}, the opposite inequality is always true,
condition \eqref{QUM} is satisfied, i.e., $\bx$ is a local quasiuniform minimum of $\varphi_1+\varphi_2$ on $B_\de(\bx)$.

Finally, let $\bar x$ be a stationary point of $\varphi_1+\varphi_2$.
By \cref{lem:characterization_stationarity},
for any $\varepsilon>0$, there exists a $\de_\varepsilon>0$ such that, for any $\de\in(0,\de_\varepsilon)$,
$\bar x$ is an $\varepsilon\delta$-minimum of $\varphi_1+\varphi_2$ on $B_\de(\bar x)$.
By assumption, there exists a $\bar\de>0$ such that $(\varphi_1,\varphi_2)$ is
quasiuniformly lower semicontinuous on $B_{\de'}(\bar x)$ for any $\de'\in(0,\bar\de)$.
For any $\varepsilon>0$, we set $\bar\de_\varepsilon:=\min(\de_\varepsilon,\bar\de)$.
Then, for any $\de\in(0,\bar\de_\varepsilon)$, as shown above,
$\bx$ is a quasiuniform $\eps\de$-minimum of $\varphi_1+\varphi_2$ on $B_\de(\bx)$, i.e., $\bx$ is a quasiuniform stationary point of $\varphi_1+\varphi_2$.
\end{proof}

\begin{remark}
In view of \cref{P3.1}\,\ref{P3.1.4}, \cref{C4.3}, \cref{P4.5}\,\ref{P4.5.3} and \cref{P4.6},
assertion \ref{P6.3.1} of \cref{P6.3} strengthens
\cite[Proposition~3.3.2]{BorZhu05}, \cite[Proposition~3.2]{Iof12} and  \cite[Proposition~2.3]{Las01}.
\end{remark}

Next, we deal with the derivation of primal necessary
conditions characterizing a quasiuniform  $\varepsilon$-minimum of the sum of two functions.

\begin{theorem}
\label{T6.1}
Let $X$ be complete,
$\varphi_1,\varphi_2$ lower semicontinuous,
$\eps>0$ and $\de>0$.
Suppose that
$\bx$ is a quasiuniform $\eps$-minimum
of $\varphi_1+\varphi_2$ on $B_{\de}(\bx)$.
Further, let
$\varphi_1$ and $\varphi_2$ be
bounded from below on ${B}_{\de}(\bx)$.
Then, for any
{sufficiently large $\rho\in(0,\de)$}
and any $\eta>0$,
there exist a number $\ga>0$ and points $\hat x_1,\hat x_2\in X$
such that
\begin{subequations}
\label{T6.1-0}
\begin{align}
\label{T6.1-1}
d((\hat x_1,\hat x_2),(\bx,\bx))<\rho,\quad
&d(\hat x_1,\hat x_2)<\eta,\quad
\varphi_\ga(\hat x_1,\hat x_2)\le(\varphi_1+\varphi_2)(\bx),
\\
\label{T6.1-2}
\sup_{\substack{u_1,u_2\in\overline{B}_{\rho}(\bar x)\\
(u_1,u_2)\ne(\hat x_1,\hat x_2)}}
&\frac{\varphi_\ga(\hat x_1,\hat x_2)-\varphi_\ga(u_1,u_2)}
{d((u_1,u_2),(\hat x_1,\hat x_2))}
<\frac{2\eps}\de,
\end{align}
\end{subequations}
where
\begin{align}
\label{phiga}
\forall u_1,u_2\in X\colon\qquad
\varphi_\ga(u_1,u_2)&:= \varphi_1(u_1)+\varphi_2(u_2)+\ga\,d(u_1,u_2).
\end{align}
\end{theorem}

\begin{proof}
Since $\bx$ is a quasiuniform $\eps$-minimum of $\varphi_1+\varphi_2$ on $B_{\de}(\bx)$,
it is a quasiuniform $\eps'$-mi\-nimum of $\varphi_1+\varphi_2$ on $B_{\de}(\bx)$ for some $\eps'\in(0,\eps)$.
Choose any numbers $\rho\in(\de\eps'/\eps,\de)$ and $\eta>0$.
Thus, $\eps/\de>\eps'/\rho$.
Set
\begin{align}
\label{T6.1P01}
\al:=\eps'/\rho^2\AND
\xi:=2(\eps/\de-\eps'/\rho).
\end{align}
By the boundedness assumption,
there is a number $c>0$ such that
\begin{align}
\label{T6.1P02}
&\forall u_1,u_2\in {B}_{\de}(\bx)\colon\quad  \varphi_1(u_1)+\varphi_2(u_2)>(\varphi_1+\varphi_2)(\bx)-c.
\end{align}
By \cref{D6.1}\,\ref{D6.1.3}, \eqref{M} and assertions \ref{L2.1.5} and \ref{L2.1.8} of \cref{L2.1}, there is a number $\ga>c/\eta$ such that
\begin{equation}
\label{T6.1P03}
\begin{aligned}
\forall u_1,u_2\in\overline B_\rho(\bar x)\colon\quad
d(u_1,u_2)<c/\gamma\;\;\Rightarrow\;\;
(\varphi_1+\varphi_2)(\bx)< \varphi_1(u_1)+\varphi_2(u_2)+\eps'.
\end{aligned}
\end{equation}
Consider a function
${\widehat\varphi_{\gamma}}\colon X\times X\to\R_\infty$
defined by
\begin{gather}
\label{T6.1P04}
\forall u_1,u_2\in X\colon\qquad
\widehat\varphi_{\gamma}(u_1,u_2):=
\varphi_\ga(u_1,u_2)+\al d((u_1,u_2),(\bx,\bx))^2.
\end{gather}
Observe that $\widehat\varphi_{\gamma}(\bx,\bx)=\varphi_{\gamma}(\bx,\bx) =(\varphi_1+\varphi_2)(\bx)$,
and $\widehat\varphi_{\gamma}$ is bounded from below on
$\overline{B}_{\rho}(\bar x)\times\overline{B}_{\rho}(\bar x)$ thanks
to \eqref{T6.1P02}.
Noting that $X\times X$ is a complete metric space,
by Ekeland variational principle (see \cref{lem:Ekeland}),
applied to the restriction of $\widehat\varphi_{\ga}$ to the closed set
$\overline{B}_{\rho}(\bar x)\times\overline{B}_{\rho}(\bar x)$,
there exist points
$\hat x_1,\hat x_2\in\overline{B}_{\rho}(\bx)$ such that
\begin{subequations}
\label{T6.1P05}
\begin{align}
\label{T6.1P05-1}
&\widehat\varphi_{\ga}(\hat x_1,\hat x_2)
\le (\varphi_1+\varphi_2)(\bx),\\
\label{T6.1P05-2}
\forall u_1,u_2\in\overline{B}_{\rho}(\bar x)\colon\quad
&\widehat\varphi_{\ga}(u_1,u_2)+\xi d((u_1,u_2),(\hat x_1,\hat x_2))
\ge
\widehat\varphi_{\ga}(\hat x_1,\hat x_2).
\end{align}
\end{subequations}
In view of \eqref{T6.1P04}
condition \eqref{T6.1P05-1} yields the last estimate in \eqref{T6.1-1}.
Moreover, it follows from \eqref{phiga}, \eqref{T6.1P02}, \eqref{T6.1P04} and \eqref{T6.1P05-1} that
\begin{equation}
\label{T6.1P07}
	\ga\,d(\hat x_1,\hat x_2)+\al d((\hat x_1,\hat x_2),(\bx,\bx))^2
	\le
	(\varphi_1+\varphi_2)(\bx)-\varphi_1(\hat x_1)-\varphi_2(\hat x_2)
	<
	c.
\end{equation}
Hence, $d(\hat x_1,\hat x_2)<c/\ga<\eta$ yielding the second estimate in \eqref{T6.1-1}.
Moreover,
by \eqref{T6.1P03} and \eqref{T6.1P07},
\begin{align*}
	\al d((\hat x_1,\hat x_2),(\bx,\bx))^2
	\le
	(\varphi_1+\varphi_2)(\bx)-\varphi_1(\hat x_1)-\varphi_2(\hat x_2)
	<
{\eps'=}
	\al\rho^2.
\end{align*}
Hence, $d((\hat x_1,\hat x_2),(\bx,\bx))<\rho$, i.e., the first estimate in \eqref{T6.1-1} holds true.
In view of \eqref{T6.1P04} we have
for any $u_1,u_2\in X$:
\begin{align*}
\varphi_\ga(\hat x_1,\hat x_2)&-\varphi_\ga(u_1,u_2)
\\&=
\widehat\varphi_{\ga}(\hat x_1,\hat x_2)-\widehat\varphi_{\ga}(u_1,u_2)
+
\al\bigl(d((u_1,u_2),(\bx,\bx))^2-d((\hat x_1,\hat x_2),(\bx,\bx))^2\bigr)
\\
&\le
\widehat\varphi_{\ga}(\hat x_1,\hat x_2)-\widehat\varphi_{\ga}(u_1,u_2)
\\&\qquad\qquad
+\al d((u_1,u_2),(\hat x_1,\hat x_2))
\big(d((u_1,u_2),(\bx,\bx))+d((\hat x_1,\hat x_2),(\bx,\bx))\big),
\end{align*}
and consequently, thanks to \eqref{T6.1P01} and \eqref{T6.1P05-2},
\begin{align*}
\sup_{\substack{u_1,u_2\in\overline{B}_{\rho}(\bar x)\\
(u_1,u_2)\ne(\hat x_1,\hat x_2)}}
\frac{\varphi_\ga(\hat x_1,\hat x_2)-\varphi_\ga(u_1,u_2)}
{d((u_1,u_2),(\hat x_1,\hat x_2))}
<\xi+2\al\rho=
2\left(\frac\eps\de-\frac{\eps'}\rho\right)+ \frac{2\eps'}{\rho}=
\frac{2\eps}\de,
\end{align*}
i.e., we arrive at \eqref{T6.1-2}.
\end{proof}

The next theorem presents dual (subdifferential) necessary conditions for
a quasiuniform $\eps$-minimum of the sum of two functions.
It is a consequence of \cref{T6.1}.

\begin{theorem}\label{T6.2}
Let $X$ be a Banach space,
$\varphi_1,\varphi_2$ lower semicontinuous, $\eps>0$ and $\de>0$.
Suppose that $\bx$
is a quasiuniform $\eps$-minimum of $\varphi_1+\varphi_2$ on $B_{\de}(\bx)$.
Further, let $\varphi_1$ and $\varphi_2$ be
bounded from below on ${B}_{\de}(\bx)$.
Then, for any $\eta>0$,
there exist points
$x_1,x_2\in B_\de(\bx)$
such that
\sloppy
\begin{subequations}
\label{eq:dual_necessary_conditions_clarke}
\begin{gather}
\label{T6.2-01}
\norm{x_1-x_2}<\eta,
\\
\label{T6.2-02}
\varphi_1(x_1)+\varphi_2(x_2)\le (\varphi_1+\varphi_2)(\bx),
\\
\label{T6.2-03}
\dist(0,{\sdc}\varphi_1(x_1)+{\sdc}\varphi_2(x_2))< {2\eps}/\de.
\end{gather}
\end{subequations}
If $X$ is Asplund, then, for any $\eta>0$,
there exist points
$x_1,x_2\in B_\de(\bx)$ satisfying \eqref{T6.2-01} and
\begin{subequations}\label{eq:dual_necessary_conditions_frechet}
\begin{gather}
\label{T6.2-04}
\varphi_1(x_1)+\varphi_2(x_2)< (\varphi_1+\varphi_2)(\bx)+\eta,
\\
\label{T6.2-05}
\dist(0,{\sdf}\varphi_1(x_1)+{\sdf}\varphi_2(x_2))< {2\eps}/\de.
\end{gather}
\end{subequations}
\end{theorem}

\begin{proof}
Let a number $\eta>0$ be given.
By \cref{T6.1}, for any sufficiently large $\rho\in(0,\de)$,
there exist a number $\ga>0$ and points $\hat x_1,\hat x_2\in X$ satisfying conditions \eqref{T6.1-0}, where the function
$\varphi_{\gamma}\colon X\to\R_\infty$
is given by \eqref{phiga}.
By \eqref{T6.1-1}, $(\hat x_1,\hat x_2)$ is an interior point of $\overline B_{\rho}(\bx)\times\overline B_\rho(\bx)$,
and it follows from \eqref{T6.1-2} that
\begin{align*}
\limsup_{\substack{(u_1,u_2)\to(\hat x_1,\hat x_2),\\(u_1,u_2)\neq(\hat x_1,\hat x_2)}}
\frac{\varphi_\ga(\hat x_1,\hat x_2)-\varphi_\ga(u_1,u_2)}
{\|(u_1,u_2)-(\hat x_1,\hat x_2)\|}
<\frac{2\eps}\de,
\end{align*}
and consequently, there is a number
$\hat\eps\in(0,{2\eps}/\de)$ such that
\begin{align*}
\liminf\limits_{\substack{(u_1,u_2)\to(\hat x_1,\hat x_2),\\(u_1,u_2)\neq(\hat x_1,\hat x_2)}}
\frac{\varphi_\gamma(u_1,u_2)-\varphi_{\gamma}(\hat x_1,\hat x_2)+ \hat\eps\|(u_1,u_2)-(\hat x_1,\hat x_2)\|}
{\|(u_1,u_2)-(\hat x_1,\hat x_2)\|}>0.
\end{align*}
By definition of the Fr\'{e}chet subdifferential,
the above inequality yields
that $0$ belongs to the subdifferential at $(\hat x_1,\hat x_2)$ of the function $(u_1,u_2)\mapsto\varphi_\ga(u_1,u_2) +\hat\eps\|(u_1,u_2)-(\hat x_1,\hat x_2)\|$,
or equivalently, in view of \eqref{phiga},
\begin{align}
\label{T6.2P03}
0\in{\sdf}\left(\varphi+g+h\right)(\hat x_1,\hat x_2),
\end{align}
where the functions
$\varphi,g,h\colon X\times X\to\R$ are given by
\begin{align}
\label{phi}
\forall u_1,u_2\in X\colon\quad
\varphi(u_1,u_2)&:=\varphi_1(u_1)+\varphi_2(u_2),
\\
\notag
g(u_1,u_2)&:=\ga \|u_1-u_2\|,\;\;
h(u_1,u_2):=\hat\eps\|(u_1,u_2)-(\hat x_1,\hat x_2)\|.
\end{align}
The next step is to apply to \eqref{T6.2P03} a subdifferential sum rule.
Note that $g$ and $h$ are convex and Lipschitz continuous, and, for all $u_1,u_2\in X$, the following relations hold for the respective subdifferentials:
\begin{subequations}
\begin{gather}
\label{T6.2P05}
\sd\varphi(u_1,u_2)=\sd\varphi_1(u_1)\times\sd\varphi_2(u_2),\\
\label{T6.2P06}
\sdc\varphi(u_1,u_2)=\sdc\varphi_1(u_1)\times\sdc\varphi_2(u_2),
\\
\label{T6.2P07}
(u_1^*,u_2^*)\in\sd g(u_1,u_2)
\quad\Rightarrow\quad
u_1^*+u_2^*=0,
\\
\label{T6.2P08}
(u_1^*,u_2^*)\in\sd h(u_1,u_2)
\quad\Rightarrow\quad
\|(u_1^*,u_2^*)\|\le\hat\eps.
\end{gather}
\end{subequations}

Inclusion \eqref{T6.2P03} obviously yields $0\in{\sdc}\left(\varphi+g+h\right)(\hat x_1,\hat x_2)$.
By the Clarke sum rule (see \cref{lem:SR}\,\ref{lem:SR.Clarke})
as well as \eqref{T6.2P07} and \eqref{T6.2P08}
there exists a subgradient $(x_{1}^*,x_{2}^*)\in{\sdc}\varphi(\hat x_{1},\hat x_2)$ satisfying
$\norm{x_{1}^*+x_{2}^*}\le\hat\eps<{2\eps}/\de$.
Set $x_1:=\hat x_1$ and $x_2:=\hat x_2$.
In view of \eqref{T6.1-1} and \eqref{T6.2P06} we have
$x_1,x_2\in B_\de(\bx)$, and conditions \eqref{eq:dual_necessary_conditions_clarke} are satisfied.

Suppose now that $X$ is an Asplund space, and set $\xi:=2{\eps}/\de-{\hat\eps}>0$.
By the fuzzy sum rule combined with the convex sum rule (see \cref{lem:SR}),
applied to \eqref{T6.2P03},
there exist a point
$(x_1,x_2)$
arbitrarily close to $(\hat x_1,\hat x_2)$
with $\varphi(x_1,x_2)$ arbitrarily close to $\varphi(\hat x_1,\hat x_2)$
and a subgradient $(x_{1}^*,x_{2}^*)\in{\sdf}\varphi(x_{1},x_2)$
such that, taking into account \eqref{T6.1-1}, \eqref{T6.2P07} and \eqref{T6.2P08},
the following estimates hold true:
\begin{gather*}
\|(x_1,x_2)-(\bx,\bx)\|<\rho,\quad
\|x_1-x_2\|<\eta,\quad
\varphi(x_1,x_2)<(\varphi_1+\varphi_2)(\bx)+\eta,
\\
\norm{x_{1}^*+x_{2}^*}<\hat\eps+{\xi}={2\eps}/\de.
\end{gather*}
Hence, $x_1,x_2\in B_\de(\bx)$
and, in view of \eqref{phi} and \eqref{T6.2P05}, conditions \eqref{T6.2-01} and \eqref{eq:dual_necessary_conditions_frechet} are satisfied.
\end{proof}

\begin{remark}
\begin{enumerate}
\item
Since the functions $\varphi_1$ and $\varphi_2$ in \cref{T6.1,T6.2} are assumed to be lower semicontinuous,
they are automatically bounded from below on \emph{some} neighbourhood of $\bx$.
We emphasize that \cref{T6.1,T6.2} require
$\bx$ to be a quasiuniform $\eps$-minimum of $\varphi_1+\varphi_2$, and
$\varphi_1$ and $\varphi_2$ to be
bounded from below on the \emph{same} fixed
neighbourhood of $\bx$.
\item
\cref{T6.2} generalizes and strengthens \cite[Theorem~4.5]{KruMeh22}.
\item
In view of \cref{P6.3}\,\ref{P6.3.2} the conclusions of \cref{T6.1,T6.2} are valid for the conventional $\eps$-minimum
if the functions are quasiuniformly lower semicontinuous on $B_\de(\bx)$.
\item
The proof of the first (general Banach space) part of \cref{T6.2} uses the Clarke subdifferential sum rule (\cref{lem:SR}\,\ref{lem:SR.Clarke}).
Clarke subdifferentials can be replaced in \cref{T6.2} by any subdifferentials possessing such an exact (see \cref{R2.1}\,\ref{R2.1.3}) sum rule in general Banach spaces.
One can use for that purpose, e.g., the $G$-subdifferentials of Ioffe; see \cite[Theorem~4.69]{Iof17}.
\end{enumerate}
\end{remark}

As consequences of \cref{T6.1,T6.2} we obtain primal and dual necessary conditions for a local quasiuniform
stationary point of a sum of functions.

\begin{corollary}
\label{C6.1}
Let $X$ be complete  and
$\varphi_1,\varphi_2$ be lower semicontinuous.
Suppose that $\bx$ is a quasiuniform stationary point of $\varphi_1+\varphi_2$.
Then, for any $\eps>0$, there
is a $\rho\in(0,\eps)$ such that,
for any $\eta>0$,
there exist a number $\ga>0$ and points $\hat x_1,\hat x_2\in X$
such that conditions \eqref{T6.1-1} are satisfied, and
\begin{gather*}
\sup_{\substack{u_1,u_2\in\overline{B}_{\rho}(\bar x)\\
(u_1,u_2)\ne(\hat x_1,\hat x_2)}}
\frac{\varphi_\ga(\hat x_1,\hat x_2)-\varphi_\ga(u_1,u_2)}
{d((u_1,u_2),(\hat x_1,\hat x_2))}
<\eps,
\end{gather*}
where the function $\varphi_\ga\colon X\to\R$ is defined by \eqref{phiga}.
\end{corollary}

\begin{corollary}
\label{C6.2}
Let $X$ be a Banach space  and
$\varphi_1,\varphi_2$ be lower semicontinuous.
Suppose that $\bx$ is a quasiuniform stationary point of $\varphi_1+\varphi_2$.
Then, for any $\varepsilon>0$, there exist points $x_1,x_2\in X$ such that
conditions \eqref{SR-1} and \eqref{T6.2-02} are satisfied, and
	\begin{gather}
\label{C6.2-2}
\dist(0,{\sdc}\varphi_1(x_1)+{\sdc}\varphi_2(x_2))<\eps.
	\end{gather}
If $X$ is Asplund, then, for any $\eps>0$,
there exist points
$x_1,x_2\in X$
such that conditions \eqref{SR-1} are
satisfied, and
\begin{gather}
\label{C6.2-3}
\dist(0,{\sdf}\varphi_1(x_1)+{\sdf}\varphi_2(x_2))<\eps.
\end{gather}
\end{corollary}

Below we provide a combined proof of the two corollaries.

\begin{proof}
[\it{Proof of \cref{C6.1,C6.2}}.]
Let $\eps>0$ and $\eta:=\eps/2$.
By the assumptions,
there exists a $\de\in(0,\eps)$ such that
\begin{align}
\label{C6.2P1}
\forall x\in B_\delta(\bar x)\colon\quad
&\varphi_1(x)-\varphi_1(\bx)>-\eta,\quad
\varphi_2(x)-\varphi_2(\bx)>-\eta.
\end{align}
and $\bx$ is a quasiuniform $\eta\de$-minimum of $\varphi_1+\varphi_2$ on $B_{\de}(\bx)$.
Thus, all the assumptions of \cref{T6.1,T6.2} are satisfied
with $\eps':=\eta\de$ in place of $\eps$.
Observe that $2\eps'/\de=\eps$, and almost all the conclusions follow immediately.
We only need to show that,
in the case of \cref{C6.2}, $\varphi_i(x_i)-\varphi_i(\bx)<\eps$, $i=1,2$.
Comparing conditions \eqref{T6.2-02} and \eqref{T6.2-04} in \cref{T6.2}, we see that condition \eqref{T6.2-04} is valid in the general as well as in the Asplund space setting.
By \eqref{T6.2-04} and \eqref{C6.2P1} we have
$$\varphi_1(x_1)-\varphi_1(\bx)< \varphi_2(\bx)-\varphi_2(x_2)+\eta <2\eta=\eps,$$
and similarly, $\varphi_2(x_2)-\varphi_2(\bx)<\eps$.
\end{proof}

\begin{remark}
\begin{enumerate}
\item
\label{R6.3.1}
Thanks to \cref{prop:quasiuniform_min_vs_stat},
\cref{C6.1,C6.2} are, particularly, applicable when the reference point is a local quasiuniform minimum.
\item
The dual necessary
conditions in \cref{C6.2} hold
not necessarily at the reference point,
but at some points arbitrarily close to it.
That is why such conditions are referred to as \emph{approximate} or \emph{fuzzy}.
Such
conditions hold
under very mild assumptions and also possess several interesting algorithmic applications;
see  e.g.\
\cite{
KruMor80,KruMor80.2,Ioffe83,Ioffe84,Kru85.1,Fab89, BorweinIoffe1996,BorweinZhu96,Iof00, Las01,Kru03,BorZhu05,Mor06.1,Iof17,
AndreaniMartinezSvaiter2010, AndreaniHaeserMartinez2011, BoergensKanzowMehlitzWachsmuth2019, KruMeh22,Mehlitz2020b,DeMarchiJiaKanzowMehlitz2022}
and the references therein.
\item
Condition \eqref{C6.2-3} in \cref{C6.2} represents a rather standard Asplund space approximate multiplier rule
(see e.g. \cite{KruMor80,Kru03,Mor06.1}),
while the general Banach space approximate multiplier rule \eqref{C6.2-2} in terms of Clarke subdifferentials is less common.
In fact, we do not know if it has been explicitly formulated in the literature.
Note that \cref{C6.2} does not assume one of the functions to be locally Lipschitz continuous (or even uniformly continuous)
as is common for multiplier rules in nonsmooth settings.
\item
The multiplier rules in \cref{C6.2} are deduced for a quasiuniform stationary point/local minimum, see also item \ref{R6.3.1}.
Thanks to \cref{P6.3}\,\ref{P6.3.1} and \ref{P6.3.3}, they apply to conventional stationary points/local minima when the pair of functions is quasiuniformly lower semicontinuous near the reference point.
Several sufficient conditions ensuring this property are given in \cref{P4.5,P4.6,P4.7,C4.3}.
In particular, the property holds if one of the functions is uniformly continuous (particularly if it is Lipschitz continuous) near the reference point.
With this in mind, the second part of \cref{C6.2} generalizes the conventional Asplund space approximate multiplier rule and makes it applicable in more general situations.
\item
Similar approximate multiplier rules in \cite[Proposition~4]{BorweinIoffe1996},
\cite[Theorem~2.9]{BorweinZhu96} and \cite[Theorem~3.3.1]{BorZhu05}
are established
(in appropriate $\be$-smooth or \Fr\ smooth spaces)
under stronger assumptions of uniform lower semicontinuity or firm uniform lower semicontinuity.
\end{enumerate}
\end{remark}

Let us revisit the setting in \cref{E3.1}.

\begin{example}
Let functions $\varphi_1,\varphi_2\colon\R^2\to\R_\infty$ be given by \eqref{E3.1-1}.
The point $(0,0)$ is obviously a minimum of $\varphi_1+\varphi_2$ and $\varphi_1(0,0)=\varphi_2(0,0)=0$.
As shown in \cref{E3.1}, ${\Lambda}_{X}^\dag(\varphi_1,\varphi_2)=0$.
Thus, $(0,0)$ is a quasiuniform minimum of $\varphi_1+\varphi_2$, and the conditions of \cref{C6.2} are satisfied.
We now check the conclusions of \cref{C6.2}.
For any $x\in\R$, we have
$$\sd\varphi_1(x,x^2)=\{(2\al x-1,-\al)\mid\al\ge0\},\quad
\sd\varphi_2(x,0)=\{(0,\al)\mid\al\ge0\}.$$
In particular,
$\sd\varphi_1(0,0)+\sd\varphi_2(0,0)=\{-1\}\times\R,$
i.e., $(0,0)\notin\sd\varphi_1(0,0)+\sd\varphi_2(0,0)$.
At the same time, if $x_1>0$, we can take $\al:=1/(2x_1)$.
Then, for any $x_2\in\R$, we have
\begin{gather*}
\varphi_1(x_1,x_1^2)-\varphi_1(0,0)=-x_1, \quad \varphi_2(x_2,0)-\varphi_2(0,0)=0,
\\
\varphi_1(x_1,x_1^2)+\varphi_2(x_2,0)=-x_1<0= \varphi_1(0,0)+\varphi_2(0,0),
\\
(0,0)=(2\al x_1-1,-\al)+(0,\al)\in \sd\varphi_1(x_1,x_1^2)+\sd\varphi_2(x_2,0).
\end{gather*}
The points $(x_1,x_1^2)$ and $(x_2,0)$ can be made arbitrarily close to $(0,0)$.
Thus, the conclusions of \cref{C6.2} are satisfied.
Observe that, when $x_1$ approaches $0$, the values $\al=1/(2x_1)$ become arbitrarily large, i.e., we are dealing with unbounded subgradients.
This is the price to pay when dropping the conventional local Lipschitz continuity assumption.

\end{example}

\section{Quasiuniform lower semicontinuity and subdifferential calculus}
\label{S7}

In this section we illustrate the value of the quasiuniform lower semicontinuity in the context
of subdifferential calculus. Our first result presents a generalized version of the
fuzzy sum rule for Fr\'{e}chet subdifferentials.

\begin{theorem}
\label{T7.1}
Let $X$ be an Asplund space,
$\varphi_1,\varphi_2\colon X\to\R_\infty$ lower semicontinuous and $\bx\in\dom\varphi_1\cap\dom\varphi_2$.
Suppose that one of the following conditions is satisfied:
\begin{enumerate}
\item\label{item:QULS}
the pair $(\varphi_1,\varphi_2)$ is firmly quasiuniformly lower semicontinuous near~$\bx$;
\item\label{item:WSLS}
$X$ is reflexive and $\varphi_1$ and $\varphi_2$ are weakly sequentially lower semicontinuous;
\item\label{item:fin_dim}
$\dim X<+\infty$.
\end{enumerate}
Then, for any $x^*\in{\sdf}(\varphi_1+\varphi_2)(\bar x)$ and $\eps>0$,
there exist points
$x_1,x_2\in X$ satisfying \eqref{SR}.
\end{theorem}

\begin{proof}
Let $x^*\in\sd(\varphi_1+\varphi_2)(\bx)$ and $\varepsilon>0$.
By definition of the \Fr\ subdifferential,
$\bx$ is a stationary point of $\varphi_1+\widehat\varphi_2$, where
$\widehat\varphi_2(x):=\varphi_2(x)-\ang{x^*,x}$ for all $x\in X$.
We now argue that the pair $(\varphi_1,\widehat\varphi_2)$ is quasiuniformly lower semicontinuous near $\bx$.
This follows from \cref{P4.4} under condition \ref{item:QULS}, and from \cref{C4.3} under condition \ref{item:WSLS}
as $\widehat\varphi_2$ is weakly sequentially lower semicontinuous if $\varphi_2$ possesses this property.
Condition \ref{item:fin_dim} is obviously a particular case of \ref{item:WSLS}.
Now, due to \cref{P6.3}\,\ref{P6.3.3}, $\bx$ is a quasiuniform stationary point of $\varphi_1+\widehat\varphi_2$.
Set $\eps':=\eps/(1+\|x^*\|)$.
By \cref{C6.2} there exist points $x_1,x_2\in B_{\eps'}(\bx)$ such that
$|\varphi_1(x_1)-\varphi_1(\bx)|<\eps'$, $|\widehat\varphi_2(x_2)-\widehat\varphi_2(\bx)|<\eps'$ and
\begin{equation}
\label{T7.1P4}
\dist(0,{\sd}\varphi_1(x_1)+{\sd}\widehat\varphi_2(x_2))<\eps'.
\end{equation}
Thus, $x_1,x_2\in B_{\eps}(\bx)$,
$|\varphi_1(x_1)-\varphi_1(\bx)|<\eps$ and
\[
	|\varphi_2(x_2)-\varphi_2(\bx)|\le |\widehat\varphi_2(x_2)-\widehat\varphi_2(\bx)|+\|x^*\|\|x_2-\bx\|< \eps'(1+\|x^*\|)=\eps.
\]
Since
${\sdf}\widehat\varphi_2(x_2)={\sdf}\varphi_2(x_2)-x^*$ (see \cref{lem:SR}\,\ref{lem:SR.0}), condition \eqref{T7.1P4} implies \eqref{SR-2}.
\end{proof}

\begin{remark}
\begin{enumerate}
\item
\cref{T7.1} strengthens \cref{lem:SR}\,\ref{lem:SR.2}.
Thanks to \cref{P4.5}\,\ref{P4.5.3}, $(\varphi_1,\varphi_2)$ is firmly quasiuniformly lower semicontinuous near $\bx$
provided that one of the functions is uniformly continuous near $\bar x$;
thus, \cref{T7.1} with condition \ref{item:QULS} also strengthens \cite[Corollary~3.4\,(ii)]{CutFab16}.

\item
An approximate sum rule, similar to \cref{T7.1} with condition \ref{item:QULS},
is established in \cite[Theorem~3.3.19]{BorZhu05}
under the stronger assumption of firm uniform lower semicontinuity (in a \Fr\ smooth space).
\end{enumerate}
\end{remark}

The next immediate corollary of \cref{T7.1} gives a
sufficient condition for the fuzzy intersection
rule for Fr\'{e}chet normals in reflexive Banach spaces.
It employs no qualification conditions and
improves the assertion of \cite[Lemma~3.1]{Mor06.1}:
it shows that 
in a reflexive space
one can always take $\la:=1$ in that lemma.

\begin{corollary}\label{cor:fuzzy_intersection_rule}
Let $X$ be a reflexive Banach space, $\Omega_1,\Omega_2\subset X$ be weakly sequentially closed and $\bar x\in\Omega_1\cap\Omega_2$.
Then, for any
$x^*\in N_{\Omega_1\cap\Omega_2}(\bar x)$ and $\eps>0$, there exist points $x_1\in \Omega_1\cap B_\eps(\bar x)$ and $x_2\in\Omega_2\cap B_\eps(\bar x)$
such that
	\[
		\dist(x^*,N_{\Omega_1}(x_1)+N_{\Omega_2}(x_2))<\eps.
	\]
\end{corollary}

\begin{remark}\label{rem:fuzzy_intersection_rule_Asplund}
\begin{enumerate}
\item
The assumptions in \cref{cor:fuzzy_intersection_rule} are trivially satisfied
for any pair of closed sets in a finite-dimensional Banach space.

\item
Employing \cref{T7.1} with condition \ref{item:QULS}, one can formulate a fuzzy intersection rule for Fr\'{e}chet normals
in an Asplund space whenever the involved sets $\Omega_1,\Omega_2$ are closed and the pair $(i_{\Omega_1},i_{\Omega_2})$ is
firmly quasiuniformly lower semicontinuous near the reference point $\bar x\in\Omega_1\cap\Omega_2$.
The latter assumption then serves as a qualification condition.
\end{enumerate}
\end{remark}

We now reinspect \cref{E5.3}.

\begin{example}
	Let $\Omega_1,\Omega_2\subset\R^2$ be given as in \cref{E5.3}.
Thus, $\Omega_1\cap\Omega_2=\{(0,0)\}$, and consequently,
$N_{\Omega_1\cap\Omega_2}(0,0)=\R^2$.
Consider the normal vector
$(1,0)\in N_{\Omega_1\cap\Omega_2}(0,0)$.
For each $x\in\R$, simple calculations show that
	\[
		N_{\Omega_1}(x,x^2)
		=
		\{(2\alpha x,-\alpha)\mid \alpha\geq 0\},
		\qquad
		N_{\Omega_2}(x,0)
		=
		\{(0,\alpha)\mid \alpha\geq 0\}.
	\]
In particular, $N_{\Omega_1}(0,0)+N_{\Omega_2}(0,0)=\{0\}\times\R$.
Hence, $(1,0)\notin N_{\Omega_1}(0,0)+N_{\Omega_2}(0,0)$.
	However, for each $x>0$, we can take $\alpha:=1/(2x)$ and find
	\[
		(1,0)
		=
		(2\alpha x,-\alpha)+(0,\alpha)
		\in
		N_{\Omega_1}(x,x^2)+N_{\Omega_2}(x,0).
	\]
When $x\downarrow 0$ we have $(x,x^2)\to (0,0)$ and $(x,0)\to (0,0)$, and consequently,
the fuzzy intersection rule indeed holds in this situation.
Recall from \cref{E5.3} that $\Omega_1$ and $\Omega_2$ are not subtransversal at $(0,0)$.
Observe that the normal vectors $(2\alpha x,-\alpha)$ and $(0,\alpha)$ become arbitrarily large as $x$ approaches $0$.

\end{example}

As another consequence of \cref{T7.1} we can derive a fuzzy chain rule in a comparatively
mild setting; cf. \cite[Section~1.2]{Kru03} and
\cite[Section~3]{MordukhovichNamYen2006}.
It employs a firm relative quasiuniform lower semicontinuity qualification condition which holds trivially, e.g.,
if the involved outer function is uniformly continuous.

\begin{corollary}
	Let $X$ and $Y$ be Asplund spaces, $\varphi\colon Y\to\R_\infty$ lower
	semicontinuous, $F\colon X\to Y$ continuous, $\bar x\in X$ and
	$F(\bar x)\in\dom\varphi$.
Suppose that one of the following conditions is satisfied:
\begin{enumerate}
\item
the function $(x,y)\mapsto\varphi(y)$ is firmly quasiuniformly lower semicontinuous relative to $\gph F$ near $(\bx,F(\bx))$;
\item
$X$ and $Y$ are reflexive, $\varphi$ is weakly sequentially lower semicontinuous and $\gph F$ is weakly sequentially closed;
\item
$\dim X<+\infty$ and $\dim Y<+\infty$.
\end{enumerate}
Then, for any $x^*\in\sdf(\varphi\circ F)(\bar x)$ and
$\varepsilon>0$, there exist points
	$\hat x\in B_\varepsilon(\bar x)$ and $\hat y\in B_\varepsilon(F(\bar x))$ such that
\begin{subequations}
\begin{gather}
\label{C7.2-1a}
|\varphi(\hat y)-\varphi(F(\bar x))|<\varepsilon,\\
\label{C7.2-1b}
		x^*\in D^*F(\hat x)(\sdf\varphi(\hat y)+\eps\mathbb B^*)+\eps\mathbb B^*.
\end{gather}
\end{subequations}
\end{corollary}

\begin{proof}
Define functions $\varphi_1,\varphi_2\colon X\times Y\to\R_\infty$ by means of
	\[
		\forall (x,y)\in X\times Y\colon\quad
		\varphi_1(x,y):=\varphi(y),\qquad
		\varphi_2(x,y):=i_{\gph F}(x,y).
	\]
By definition of the Fr\'{e}chet subdifferential
	$x^*\in\sdf(\varphi\circ F)(\bar x)$ implies that
	$(x^*,0)\in\sdf(\varphi_1+\varphi_2)(\bar x,F(\bar x))$.
\cref{T7.1} gives the existence of $\hat x\in B_\varepsilon(\bar x)$ and $\hat y\in B_\varepsilon(F(\bar x))$ satisfying \eqref{C7.2-1a} and 	
	\[
		\dist\bigl(
			(x^*,0),\{0\}\times\sd\varphi(\hat y)+N_{\gph F}(\hat x,F(\hat x))
			 \bigr)
		<
		\varepsilon.
	\]
The last condition obviously implies \eqref{C7.2-1b}.
\end{proof}

\section{An application in optimal control}
\label{S8}

We revisit the setting of \cref{ex:lower_semicontinuity_relative_to_a_set}.
Let $D\subset\R^d$ be some bounded open set.
We consider a continuously differentiable mapping $S$ from $L^2(D)$ to a Hilbert space $H$,
the so-called \emph{control-to-observation operator},
which assigns to each control function $x\in L^2(D)$ an observation $S(x)\in H$.
Typically, $S$ represents the composition of the solution operator
associated with a given variational problem
(e.g.\ a partial differential equation or a variational inequality)
and some mapping which sends the output of the variational problem to the observation space $H$.
In optimal control, a function $x$ often has to be chosen such that
$S(x)$ is close to some desired object $y_\textup{d}\in H$ which can be modeled
by the minimization of the smooth term $\tfrac12\norm{S(x)-y_\textup{d}}^2$.
There are often other requirements which have to be respected in many situations.
For example, a control has to belong to a simple constraint set
$\Omega\subset L^2(D)$
or has to be \emph{sparse}, i.e., it has to vanish on large parts of the domain $D$.
Here, we take a closer look at the sparsity-promoting function
$\varphi\colon L^2(D)\to\R$ given in \eqref{eq:sparsity_promoting_function}.
Furthermore, we assume that $\Omega$ is given as in
\eqref{eq:box_constrained_set} where $x_a,x_b\in L^2(D)$ satisfy
$x_a(\omega)< 0< x_b(\omega)$ almost everywhere on $D$.
We note that $\Omega$ is closed and convex, so the various normal cones to this
set coincide with the one in the sense of convex analysis.

We investigate the optimal control problem
\begin{equation}\label{eq:OC}\tag{OC}
	\min\{
		f(x)
		+
		\varphi(x)
		\mid
		x\in \Omega
		\}
\end{equation}
where $f\colon L^2(D)\to\R$ is an arbitrary continuously differentiable function
and keep in mind that a possible choice for $f$ would be the typical target-type
function
$x\mapsto\tfrac12\norm{S(x)-y_{\textup{d}}}^2+\tfrac{\sigma}{2}\norm{x}^2$
where $\sigma\geq 0$ is a regularization parameter.
We identify the dual space of $L^2(D)$ with $L^2(D)$.
Thus, for any $x\in L^2(D)$, $f'(x)$ can be interpreted as a function from $L^2(D)$.
Problems of type \eqref{eq:OC} were already considered e.g.\ in
\cite{ItoKunisch2014,NatemeyerWachsmuth2021,Wachsmuth2019} from the viewpoint of
necessary and sufficient optimality conditions as well as numerical solution methods.

Before we can state necessary optimality conditions for this optimization problem
it has to be clarified how the subdifferentials of $\varphi$ look like.
This has been investigated in the recent paper \cite{MehlitzWachsmuth22}.
Before presenting the formulas we need to recall the concept of so-called
\emph{slowly-decreasing} functions, see \cite[Definition~2.4, Theorem~2.10]{MehlitzWachsmuth22}
as well as the discussions therein.

\begin{definition}\label{def:SD_functions}
	A function $x\in L^2(D)$ is called \emph{slowly decreasing} if
	\[
		\lim\limits_{t\downarrow 0}
			\blambda(\{0<|x|\leq t\})/t^2
		=
		0.
	\]
\end{definition}

The following lemma can be distilled from \cite[Theorems~3.5 and 3.7]{MehlitzWachsmuth22}.

\begin{lemma}\label{lem:subdifferential_of_sparsity_promoting_function}
	For $\varphi\colon L^2(D)\to\R$ from \eqref{eq:sparsity_promoting_function}
	and $x\in L^2(D)$, we have
	\begin{align*}
		\sdf\varphi(x)
		&=
		\begin{cases}
			\{x^*\in L^2(D)\,|\,
				\{x^*\neq 0\}\subset\{x=0\}
			\}
			&
			x\text{ slowly decreasing,}
			\\
			\emptyset
			&
			\text{otherwise,}
		\end{cases}
		\\
		\sdm\varphi(x)
		&=
		\{x^*\in L^2(D)\,|\,
				\{x^*\neq 0\}\subset\{x=0\}
			\}.
	\end{align*}
\end{lemma}

In our first result we present approximate stationarity conditions for
\eqref{eq:OC}.

\begin{theorem}\label{thm:approximate_conditions_OC}
	Let $\bar x\in L^2(D)$ be a local minimum of \eqref{eq:OC}.
	Then, for each $\varepsilon>0$, there exist a slowly decreasing
	function $x_1\in B_\varepsilon(\bar x)$, some $x_2\in\Omega\cap B_\varepsilon(\bar x)$
	and $x_1^*,x_2^*\in L^2(D)$ such that
	\begin{subequations}\label{eq:approxiate_conditions_OC}
		\begin{gather}
			\label{eq:approximate_conditions_OC_stat}
				\norm{f'(x_2)+x_1^*+x_2^*}<\varepsilon,
				\\
			\label{eq:approxiate_conditions_OC_function_value}
				|\varphi(x_1)-\varphi(\bar x)|
				<
				\varepsilon,
				\\
			\label{eq:approximate_conditions_OC_sparse}
				\{x_1^*\neq 0\}\subset\{x_1=0\},
				\\
			\label{eq:approximate_conditions_OC_normal_cone}
				x_2^*\geq 0\text{ a.e.\ on }\{x_2>x_a\},\quad
				 x_2^*\leq 0\text{ a.e.\ on }\{x_2<x_b\}.
		\end{gather}
	\end{subequations}
\end{theorem}
\begin{proof}
	From \cref{ex:lower_semicontinuity_relative_to_a_set} we know that
	$(\varphi,i_\Omega)$ is firmly quasiuniformly lower semicontinuous near $\bar x$,
	and due to \cref{P4.4}, this extends to $(\varphi,f+i_\Omega)$
	since $f$ is locally Lipschitz continuous and, thus, uniformly continuous
	on each sufficiently small ball around $\bar x$.
	Applying \cref{P6.3}\,\ref{P6.3.3} shows that $\bar x$ is a quasiuniform
	stationary point of $\varphi+(f+i_\Omega)$.	
	Thus, \cref{lem:SR}\,\ref{lem:SR.0} and \cref{C6.2} yield the existence of
	$x_1\in B_\varepsilon(\bar x)$ and $x_2\in \Omega\cap B_\varepsilon(\bar x)$
	satisfying \eqref{eq:approxiate_conditions_OC_function_value} and
	\[
		\dist(-f'(x_2),\sdf\varphi(x_1)+N_\Omega(x_2))<\varepsilon.
	\]
	Now, the remaining assertions of the theorem follow from
	\cref{lem:subdifferential_of_sparsity_promoting_function}
	and the well-known characterization of the normal cone
	$N_\Omega(x_2)$.
\end{proof}

We now take the limit as $\varepsilon\downarrow 0$ in the system
\eqref{eq:approxiate_conditions_OC} in order to obtain a conventional
stationarity condition.

\begin{theorem}\label{thm:sharp_conditions_OC}
	Let $\bar x\in L^2(D)$ be a local minimum of \eqref{eq:OC}.
	Then
	\begin{subequations}\label{eq:sharp_conditions_OC}
		\begin{align}
			\label{eq:sharp_conditions_OC_inactive}
				f'(\bar x) &=0	\quad\text{a.e.\ on }\{\bar x\neq 0\}\cap\{x_a<\bar x<x_b\},
				\\
			\label{eq:sharp_condition_OC_lower_bound}
				f'(\bar x) &\geq 0	\quad\text{a.e.\ on }\{\bar x=x_a\},\\
			\label{eq:sharp_condition_OC_upper_bound}
				f'(\bar x)	&\leq 0	\quad\text{a.e.\ on }\{\bar x=x_b\}.
		\end{align}
	\end{subequations}
\end{theorem}
\begin{proof}
	For each $k\in\N$, we apply \cref{thm:approximate_conditions_OC} with
	$\varepsilon:=1/k$ in order to find a slowly decreasing function
	$x_{1k}\in B_{1/k}(\bar x)$, some $x_{2k}\in \Omega\cap B_{1/k}(\bar x)$
	and $x_{1k}^*,x_{2k}^*\in L^2(D)$ such that
	\begin{subequations}\label{eq:sequential_stationarity_OC}
	\begin{gather}
		\label{eq:sequential_stationarity_OC_stat}
		\norm{f'(x_{2k})+x_{1k}^*+x_{2k}^*}<1/k,\\
		\label{eq:sequential_stationarity_OC_sparse}
		\{x_{1k}^*\neq 0\}\subset\{x_{1k}=0\},\\
		\label{eq:sequential_stationarity_normal_cone}
		x_{2k}^*\geq 0\text{ a.e.\ on }\{x_{2k}>x_a\},\quad
		 x_{2k}^*\leq 0\text{ a.e.\ on }\{x_{2k}<x_b\}.
	\end{gather}
	\end{subequations}
	Clearly, we have $x_{1k}\to\bar x$ and $x_{2k}\to\bar x$ as $k\to+\infty$.
Set $x_k^*:=x_{1k}^*+x_{2k}^*$ for each $k\in\N$.
	Due to \eqref{eq:sequential_stationarity_OC_stat}, we find $x_k^*\to-f'(\bar x)$
	from continuous differentiability of $f$.
	Along a subsequence (without relabeling) we can assume that the convergences
	$x_{1k}(\omega)\to\bar x(\omega)$, $x_{2k}(\omega)\to\bar x(\omega)$ and
	$x_k^*(\omega)\to-f'(\bar x)(\omega)$  hold for almost all $\omega\in D$.
	
	Note that \eqref{eq:sequential_stationarity_normal_cone} gives
	\[
		x_{2k}^*(\omega)
		\begin{cases}
			\leq 0	&	\omega\in\{x_{2k}=x_a\},\\
			\geq 0	&	\omega\in\{x_{2k}=x_b\},\\
			= 0		&	\omega\in\{x_a<x_{2k}<x_b\}
		\end{cases}
	\]
	for almost every $\omega\in D$.
	For almost every $\omega\in \{\bar x=x_a\}$, we have $x_{1k}(\omega)<0$ as well as
	$x_{2k}(\omega)<0$ and, thus, $x_k^*(\omega)\leq 0$ for large enough $k\in\N$, see
	\eqref{eq:sequential_stationarity_OC_sparse} as well, and taking
	the limit as $k\to +\infty$ gives \eqref{eq:sharp_condition_OC_lower_bound}.
	Similarly, we can show \eqref{eq:sharp_condition_OC_upper_bound}.
	Finally, for almost every $\omega\in\{\bar x\neq 0\}\cap\{x_a<\bar x<x_b\}$, we have
	$\omega\in\{x_{1k}\neq 0\}$ and $\omega\in\{x_a<x_{2k}<x_b\}$ for large enough $k\in\N$,
	giving $x_k^*(\omega)=0$ for any such $k\in\N$, and taking the limit
	gives \eqref{eq:sharp_conditions_OC_inactive}.
\end{proof}

\begin{remark}\label{rem:optimality_conditions_OC}
	\begin{enumerate}
		\item Due to \cref{lem:subdifferential_of_sparsity_promoting_function}, we
			can reformulate \eqref{eq:sharp_conditions_OC} in terms of multipliers.
			If $\bar x\in L^2(D)$ is a local minimum of \eqref{eq:OC}, then there
			exist $x_1^*\in\sdm\varphi(\bar x)$ and $x_2^*\in N_\Omega(\bar x)$
			such that $f'(\bar x)+x_1^*+x_2^*=0$.
		\item The set of slowly decreasing functions is not closed in $L^2(D)$.
			It is, thus, not surprising that taking the limit as $\varepsilon\downarrow 0$ in
			\eqref{eq:approxiate_conditions_OC} annihilates this information.
		\item The stationarity conditions from \eqref{eq:sharp_conditions_OC} clearly
			promote sparse solutions since only on $\{\bar x=0\}$, there is lots of
			freedom available regarding the sign of $f'(\bar x)$.
		\item Slightly more restrictive necessary optimality conditions can be obtained
			via the Pontryagin maximum principle if $f$ satisfies the additional requirement
			\[
				f(x)=f(\bar x)+f'(\bar x)(x-\bar x)+\oo\bigl(\norm{x-\bar x}_{L^1(D)}\bigr)
			\]
			for each $x\in \Omega$, see
			\cite[Theorem~2.5]{Wachsmuth2019}.
			Here, $L^1(D)$ is the space of (equivalence classes of) integrable functions
			equipped with the usual norm $\norm{\cdot}_{L^1(D)}$.
			By boundedness of $D$, $L^2(D)$ is continuously embedded in $L^1(D)$.
	\end{enumerate}
\end{remark}

\section{Conclusions}\label{sec:conclusions}

In this paper we reinspect the
popular
decoupling approach
that has proved to be useful in various areas of nonlinear analysis and optimization
involving problems which can be modeled as
minimization of the sum of two extended-real-valued functions.
Several
decoupling quantities from the literature like the \emph{uniform} or \emph{decoupled infimum} are reviewed, and
some new ones like the \emph{quasiuniform infimum} are introduced and the relations between them are studied.
We exploit these decoupling tools to define the concepts of (\emph{firm}) \emph{uniform} and \emph{quasiuniform lower semicontinuity}
describing certain stability properties of the infimum of the
sum whenever
the latter is decoupled.
The relationship between these uniform lower semicontinuity properties is studied
and sufficient conditions for their validity are established.
Our new concepts and results are embedded into the rich landscape of related literature which addresses the decoupling approach.

These abstract findings are used for the derivation of primal and dual necessary conditions
characterizing stationary points (and, in particular, local minimizers) of sums
of functions under mild assumptions.
Consequences for the calculus of Fr\'{e}chet subdifferentials and normals are distilled.
For instance, mild sufficient conditions for the validity of the fuzzy sum rule
for Fr\'{e}chet subdifferentials are given which apply in situations where the involved summands
are not even uniformly continuous around the reference point.
An illustrative example from sparse optimal control visualizes applicability of our results
within a typical infinite-dimensional setting.

\section*{Acknowledgements}
The authors acknowledge fruitful discussions with Gerd Wachsmuth which led to
the construction of \cref{E3.4}.
We wish to thank the referee for helpful comments which allowed us to improve the presentation.

Research is supported by the
Australian Research Council, project DP160100854; the DFG Grant Bilevel Optimal Control: Theory, Algorithms,
and Applications, grant WA 3636/4-2, within the Priority Program SPP 1962 (Non-smooth and Complementarity-based Distributed Parameter Systems: Simulation and Hierarchical Optimization); and
grants GA\v{C}R 20-22230L and by RVO: 67985840.
The second author benefited from the
support of the European Union's Horizon 2020 research and innovation programme under the Marie Sk{\l}odowska-Curie
Grant Agreement No. 823731 CONMECH.
\medskip

\noindent
{\bf Declarations of interest}: none

\addcontentsline{toc}{section}{References}


\end{document}